\documentclass[12pt,reqno]{amsart}
\usepackage[margin=1in]{geometry}  
\usepackage[english]{babel}
\usepackage[utf8]{inputenc} 
\usepackage{hyperref} 
\usepackage{amsmath,amssymb,amsthm,bbm}
\usepackage[makeroom]{cancel}
\usepackage{cases}
\usepackage{mathtools}
\usepackage{esint} 
\usepackage{graphicx, color}
\usepackage[dvipsnames]{xcolor} 
\usepackage{enumitem}  
\usepackage{pgf}
\newcommand{\dx}{{\rm d} x}

\newcommand{\dt}{{\rm d} t}
\newcommand{\ds}{{\rm d} s}

\renewcommand{\d}{{\rm d}}

\newcommand{\numberset}{\mathbb}
\newcommand{\R}{\numberset{R}}

\newcommand{\N}{\numberset{N}}
\newcommand{\C}{\numberset{C}}

\newcommand{\Rz}{\R}
\newcommand{\SO}{{\rm SO}}

\newcommand{\GL}{{\rm GL}}
\newcommand{\haz}{\widehat}

\newcommand{\RR}{\mathcal{R}}
\newcommand{\HH}{\mathcal{H}}
\newcommand{\W}{\mathcal{W}}

\newcommand{\DD}{\mathbb{D}}

\newcommand{\eps}{\varepsilon}

\newcommand{\J}{V^{J}}
\newcommand{\cJ}{\cW}
\newcommand{\cW}{c_W}

\DeclareMathOperator*{\argmin}{arg\,min}
\newcommand{\weak}{\rightharpoonup}

\theoremstyle{plain}
\newtheorem{theorem}{Theorem}[section]
\newtheorem{definition}{Definition}[section]
\newtheorem{proposition}{Proposition}[section]
\newtheorem{corollary}{Corollary}[section]

\theoremstyle{definition}

\numberwithin{equation}{section}
 
\title[   Viscoelasticity and accretive phase-change]{  
     Viscoelasticity and accretive phase-change\\ at finite strains}

\author[A. Chiesa] {Andrea Chiesa} 
\address[Andrea Chiesa]{University of Vienna, Faculty of Mathematics  and Vienna School of Mathematics,
                Oskar-Morgenstern-Platz 1, A-1090 Vienna, Austria}
\email{andrea.chiesa@univie.ac.at}
	\urladdr{http://www.mat.univie.ac.at/$\sim$achiesa}
\author[U. Stefanelli]{Ulisse Stefanelli} 
	\address[Ulisse Stefanelli]{University of
		Vienna, Faculty of Mathematics,
                Oskar-Morgenstern-Platz 1, A-1090 Vienna, Austria, 
		University of Vienna, Vienna Research Platform on Accelerating
		Photoreaction Discovery, W\"ahringerstra\ss e 17, 1090
                Wien, Austria, and Istituto di	Matematica Applicata e Tecnologie Informatiche {\it E. Magenes}, via
		Ferrata 1, I-27100 Pavia, Italy}
	\email{ulisse.stefanelli@univie.ac.at}
	\urladdr{http://www.mat.univie.ac.at/$\sim$stefanelli}


\begin{document}

\maketitle

\begin{abstract}
  We investigate    the evolution of   a two-phase
  viscoelastic material at finite strains.    The   phase evolution is assumed
  to be
irreversible   :    One   phase accretes in time in its normal
direction, at the expense of the
other. Mechanical response depends on the phase.    At the same
time,   growth is influenced by the mechanical state at the boundary of the
accreting phase, making the model fully coupled. This setting is
inspired by the early stage development of solid
tumors, as well as by the swelling of polymer gels. We formulate the evolution problem by coupling the
balance of momenta in weak form and the growth dynamics
in the viscosity sense. Both a diffused- and a sharp-interface variant
of the model
are proved to admit solutions    and the sharp-interface limit is
investigated.  

 \end{abstract}

\section{Introduction}
 This paper is concerned with the evolution of a
 viscoelastic compressible solid undergoing phase    change.   We assume that the
 material presents two phases, of which one grows at the expense of
 the other by {\it accretion}. In particular, the phase-transition    front  
 evolves in a {\it normal} direction to the accreting phase, with a
    growth  
 rate depending on the deformation. This behavior is indeed common to
 different material systems. It may be observed in the early stage development of
 solid tumors \cite{tumor1,tumor2,tumor3}, where the neoplastic tissue invades the
 healthy one. Swelling in polymer gels also follows a similar
 dynamics,
 with the swollen phase accreting in the dry one
 \cite{Lucantonio,Shibayama} and causing a volume increase.
 Accretive growth can be observed in some
 solidification processes~\cite{cellulardendritic,Wang
   image}, as well.
 
  The focus of the modelization is on describing the  interplay between
  mechanical deformation and accretion  .   On the one hand, the two
phases are assumed to have a different mechanical response, having an
effect on the viscoelastic evolution of the medium. On the other hand, the time-dependent mechanical
deformation is assumed to influence the growth process, as is indeed common in biomaterials \cite{Goriely}, polymeric gels \cite{Yamaue}, and
solidification \cite{solidification1}. The mechanical and phase
evolutions are thus fully coupled.


 The state of the system is described by the pair $(y,\theta) : [0,T] \times U \to \Rz^d \times
 [0,\infty)$, where $T>0$ is some final time and $U   \subset \Rz^d\,
 (   d\geq 2  )$ is the
 reference configuration of the body. Here, $y$ is the deformation
 of the medium while $\theta$ determines its phase. More precisely, for
 all $t\in [0,T]$ the {\it accreting (growing) phase} is identified as the
 sublevel $\Omega(t):=\{x \in U \ | \ \theta(x)<t\}$, whereas the {\it
    receding 
   phase} corresponds to $U \setminus \overline{\Omega(t)}$. The value 
$\theta(x)$ formally corresponds to  the
 time at which the point $x\in U$ is added to the growing phase. As
 such, $\theta$ is usually referred to as \emph{time-of-attachment}
 function. An illustration of the    notation   is given in Figure
 \ref{fig:setting}.
 \pgfdeclareimage[width=150mm]{setting}{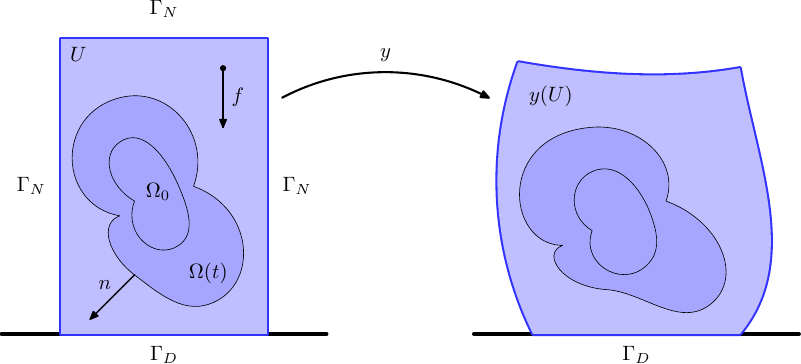}
 \begin{figure}[h]
   \centering
   \pgfuseimage{setting}
   \caption{   Illustration of the
notation in the reference domain (left) and in the
deformed one (right).}
   \label{fig:setting}
 \end{figure}

As growth processes and mechanical equilibration typically occur on very
different time scales, we    neglect inertial effects and assume the
evolution to be
viscoelastic.   
This calls for specifying the stored
 energy density $W(\theta(x){-}t,\nabla y)$ and the viscosity $R(\theta(x){-}t,\nabla y,\nabla
 \dot y)$ of the medium, as well as the applied body forces
 $f(\theta(x){-}t,x)$. All these quantities are assumed to be
 dependent on the phase via the {\it sign} of $\theta(x){-}t$, which
 indeed distinguishes the two phases, in the
 spirit of the celebrated {\it level-set method}
 \cite{Osher,Sethian}. In addition, we include a
 {\it second-gradient} regularization term in the energy of the form
 $H(\nabla^2y)$, which we take to be phase independent, for
 simplicity. All in all, the    viscoelastic evolution   system takes the
 form
 \begin{equation}
   \label{eq:y_intro}
   {}-\operatorname{div}\left(\partial_{\nabla y}W (\theta(x){-}t, \nabla
     y) +\partial_{\nabla \dot y}R(\theta(x){-}t,\nabla y,\nabla \dot y)-\operatorname{div}{\rm D}  H(\nabla^2y)\right)= f(\theta(x){-}t,x).
 \end{equation}
This system is solved weakly, complemented by mixed boundary conditions on $y$
and a homogeneous natural condition on the {\it hyperstress} $   {\rm D}
H(\nabla^2 y)$, namely, 
\begin{align}
  &y={\rm id} \ \ \text{on}  \ \ [0,T]\times\Gamma_D,\label{eq:boundary condition}\\
  &{\rm D}H(\nabla^2y) {:}  (\nu\otimes \nu)=0 \ \  \text{on}  \ \ [0,T]\times\partial U,\label{eq: boundary condition H}\\
 &\left(\partial_{\nabla y} W (\theta(x){-}t, \nabla
   y){+}\partial_{\nabla \dot y}R_\eps(\theta(x){-}t,\nabla y, \nabla \dot y)\right)\nu \notag\\
&\quad -{\rm div}_S\,({\rm
  D}H(\nabla^2 y)    \nu  )   =0  \ \ \text{on} \ \ [0,T]\times \Gamma_N,\label{eq: neumann condition}
\end{align}
 where $\nu$ is the outer unit normal to $\partial U$, $\Gamma_D$ and $
\Gamma_N$ 
are the Dirichlet and Neumann part of the boundary $\partial U$,
respectively, and ${\rm div}_S$ denotes the surface divergence on
$\partial U$    \cite{MielkeRoubicek}.   

The    viscoelastic evolution   system is coupled to the phase evolution
by requiring that the time-of-attachment function $\theta$ solves 
the generalized eikonal equation
\begin{equation}
  \label{eq:theta_intro}
  \gamma\big(y(\theta(x)   \wedge T  ,x),\nabla y(\theta(x)    \wedge T  ,x)\big)|\nabla (-\theta)(x)|=1
\end{equation}
for all $x$ in the complement of a given initial set  $\Omega_0\subset\subset
U$ where we set $\theta=0$. This corresponds to assuming that $\Omega (t)$ accretes in its normal
direction, with {\it growth rate} $\gamma(\cdot)>0$. 
 More precisely, the evolution of the generic point $x(t)$
 on the boundary
 $\partial \Omega(t) $ follows the ODE flow
 \begin{equation*}
    \frac{\d}{\dt}   x(t) = \gamma \big(y(t,x(t)), \nabla y(t,x(t))\big) \,
    \nu (x(t)) 
 \end{equation*}
 where $   \nu (x(t))$ indicates the normal to $\partial \Omega (t)$ at
 $x(t)$. Accretive growth is paramount to a wealth of different biological
 models \cite{Taber}, including plants and trees \cite{plants,trees} and the formation of
 hard tissues like horns or shells
 \cite{Moulton,calcific1,calcific2}. The dependence of the growth rate
 $\gamma$ on the actual position  and strain
 is intended to model the possible influence of local features such as
 nutrient concentrations, as well as of the local mechanical state \cite{Goriely}. 
 Note that accretive growth occurs in a variety of
 nonbiological systems, as well.  These include crystallization
 \cite{crystallization0,crystallization2}, sedimentation of rocks \cite{Ganghoffer2}, glacier
 formation, accretion of celestial bodies \cite{planet}, as well as technological
 applications, from epitaxial deposition \cite{Masi}, to coating, masonry, and 3D
 printing \cite{3Dprint1,3Dprint2}, just to mention a few. 

 By assuming smoothness and differentiating the equation $\theta(x(t))=t$
 in time one obtains $\nabla \theta(x(t)){\cdot}  \frac{\d}{\dt}    x(t)=1$. This,
 together with the above flow rule for $x(t)$  
 and $   \nu  (x(t)) = \nabla
 \theta (x(t))/|\nabla
 \theta (x(t))|$, originates the generalized eikonal equation \eqref{eq:theta_intro}.
    As    the growth rate $\gamma$ in \eqref{eq:theta_intro} depends
 on the actual deformation $y(\theta(x)   \wedge T  ,x)$ and strain $\nabla
 y(\theta(x)    \wedge T  ,x)$ at the growing interface,
 system \eqref{eq:y_intro}--\eqref{eq:theta_intro}    is   fully
 coupled.

 We specify the initial conditions for the system by setting
 \begin{align} 
  &\theta  = 0
   \ \ \text{on} \ \Omega_0,\label{eq:initial_intro0}\\
   &y(0,\cdot) = y_0 \ \ \text{on} \  U, \label{eq:initial_intro}
 \end{align}
 where the initial deformation  $y_0$ and the initial portion of
 the growing phase $\Omega_0$ are given. Note that $\Omega_0$ is not
 required to be connected, in order to possibly model the onset of the
 accreting phase at different sites. On the other hand, the evolution
 described by \eqref{eq:theta_intro} does not preserve the topology
 and disconnected accreting regions may    eventually   coalesce over time.

 The aim of this paper is to present an existence theory to the
 initial and boundary value problem
 \eqref{eq:y_intro}--\eqref{eq:initial_intro}. We tackle
 both a {\it sharp-interface} and a {\it diffused-interface} version
 of the model,
 by tuning the assumptions on $W$ and $R$, see Sections
 \ref{sec:energy}--\ref{sec:dissipation}. More precisely, in the
 diffused-interface model we assume that energy and dissipation
 densities  change smoothly as functions of the phase indicator
 $\theta(x){-}t$ across a region of width $\eps>0$, namely for
 $-\eps/2<\theta(x)-t<\eps/2$.  On the contrary, in the 
 sharp-interface case material potentials are assumed to be discontinuous 
across the phase-change surface $\{\theta(x)=t\}$.
 
In both regimes, we prove that the fully coupled system
\eqref{eq:y_intro}--\eqref{eq:initial_intro} admits a weak/viscosity
solution, see Definition \ref{def:weak sol} below. More precisely, the 
viscous evolution \eqref{eq:y_intro}--\eqref{eq: neumann condition} is solved weakly, whereas the growth
dynamics equation \eqref{eq:theta_intro} is solved in the viscosity
sense, see Theorem \ref{Thm:existence sol}. We moreover prove that solutions fulfill the energy equality,
where the energetic contribution of the phase transition is
characterized, see Proposition \ref{prop:energy}.  
As a by-product, solutions of the
diffused-interface model for $\eps>0$ are proved to converge up to subsequences to
solutions of the sharp-interface model as the parameter $\eps$
converges to $0$, see Corollary~\ref{prop:sharp}. 

 Before going on, let us mention that the engineering literature on
growth mechanics is vast. Without any claim of completeness, we
mention \cite{linearizedlit,TruskinovskyZurlo2}  and
\cite{finitestrainlit,WalkingtonDayal1,WalkingtonDayal2, Sozio,
  TruskinovskyZurlo1} as examples of linearized and finite-strain
theories, respectively. On the other hand, mathematical existence theories in growth
mechanics are scant, and we refer to
\cite{Bangert,existence2,Ganghoffer} for some recent results. To the
best of our knowledge, no existence result for finite-strain accretive-growth
mechanics is currently available. In the linearized case, an existence
result for the model \cite{TruskinovskyZurlo2} has been obtained in 
\cite{DavoliNikStefanelliTomassetti}.

The paper is structured as follows. Section \ref{sec: setting} is
devoted to the statement of the main existence result, Theorem
\ref{Thm:existence sol}. In Section \ref{sec:energy_equality}, we give
the proof of the energy identity. The proof of Theorem
\ref{Thm:existence sol} is then split in
Sections \ref{sec: Diffusive model} and \ref{Sec: eps to zero}, respectively focusing on the
diffused-interface and the sharp-interface setting. In the
diffused-interface case, the proof relies on an iterative
construction, where the mechanical and the growth problems are solved
in alternation. The existence proof for the sharp-interface model is
obtained by    taking the limit as  $\eps\to 0$ in the
diffused-interface model.

\section{Main results}\label{sec: setting}

 In this section, we specify assumptions, introduce the
weak/viscosity notion of solution, and state the main results for
problem \eqref{eq:y_intro}--\eqref{eq:initial_intro}.

\subsection{Notation}\label{sec: notation}
In what follows, we denote by $\R^{d\times d}$ the Euclidean space of
$d{\times} d$ real matrices,  $d   \geq 2 $,    by  $\R^{d\times d}_{\rm sym}$
the  subspace   of symmetric matrices, and by $I$ the identity
matrix. Given $A\in \R^{d\times d}$, we
indicate  its  transpose by $A^\top$ and  its Frobenius
  norm by $
    |A|^2 := A{:}A$,
  where the contraction product between matrices  $A,\,B\in
  \R^{d\times d}$  is defined as $A {:}
   B  :=  A_{ij}B_{ij}$  (we use the summation convention over
   repeated indices).  Analogously, let $\R^{d\times d\times d}$ be the set of real $3$-tensors, and define their contraction product as $A \smash{\vdots}
   B  := A_{ijk}B_{ijk}$ for $A,B\in \R^{d\times d\times d}$.    A 4-tensor
  $\C\in\R^{d\times d\times d\times d}$ is said to be {\it major symmetric}
  if $\C_{ijk\ell} = \C_{k\ell ij}$ and {\it minor symmetric} if
  $\C_{ijk\ell} = \C_{ij\ell k}=\C_{jik\ell}$. Given a major and minor   symmetric positive definite 4-tensor
  $\C\in\R^{d\times d\times d\times d}$ and the matrix $A\in
  \R^{d\times d}$ we indicate by $\C{:}A\in \R^{d\times d}$ and $ A{:}\C\in \R^{d\times d}$ the matrices
  given in components by $(\C{:}A)_{ij} = \C_{ijk\ell}A_{k\ell}$ and
  $(A{:}\C)_{ij} = A_{k\ell}C_{k\ell i j}$, respectively.   
  We shall use the following matrix sets $SO(d):=\{ A \in   
  \R^{d\times d} \;|\; \det{A}=1 \ \text{and}   \  A A^{\top}
   =I\}$ and $\GL_{+}(d):=\{ A  \in \R^{d\times d} \;|\; \det  A >0
   \}$.
   
   The scalar product of two vectors $a,\, b\in \R^d$ is
  classically indicated by $a{\cdot}b$. 
  The symbol $B_R\subset \R^{d}$ denotes the open ball of radius $R>0$
  and center $0\in\R^{d}$, $|E|$ indicates the Lebesgue measure
  of the Lebesgue-measurable set $E\subset \Rz^d$,  and
  $\mathbbm{1}_E$ is the corresponding characteristic function,
  namely, $\mathbbm{1}_E(x)=1$ for $x\in E$ and $\mathbbm{1}_E(x )=0$
  otherwise. For $E \subset \Rz^d$ nonempty and $x\in \Rz^d$ we define
  ${\rm dist}(x,E):=\inf_{e\in E}|x{-}e|$. We denote by
  $\mathcal{H}^{d-1}$ the $(d{-}1)$-dimensional Hausdorff
  measure    and   $x{\wedge}y\coloneqq \min\{x;y\}$ for all
  $x\,,y\in \R$.  

     In the following, we use the symbol $\dot u$ for the partial
  time derivative of the generic time-dependent function $u$, whereas $\frac{\d}{\dt}$
  stands for the total time derivative, in case $u$ depends on time
  only.

Henceforth, we indicate by  $c$ a  generic
  positive constant possibly depending on data but independent of
  the parameter $\eps$ and of the
  discretization step $\tau$. 
 Note that the value of $c$ may change even within
 the same line.  

\subsection{Admissible deformations}\label{sec:states}

Fix the final time $T>0$ and let the reference configuration $U
\subset \R^d \     (d\geq 2)$   be nonempty, open,
connected, and bounded. We assume that the
boundary $\partial U$ is Lipschitz, with $\Gamma_D, \,
\Gamma_N\subset \partial U$ disjoint and open in the topology of
$\partial U$,   $\Gamma_D \not = \emptyset$ and $\overline{\Gamma_D}  \cup  \overline{\Gamma_N} =
\partial U$, where the closure is taken in the topology of $\partial
U$. In the following, we use the short-hand notation $Q:=(0,T)\times U$ and
$\Sigma_D:=(0,T)\times \Gamma_D$.

Deformations are assumed to belong to the    affine space   
  \begin{align*} 
    W^{2,p}_{\Gamma_D}(U ;\R^d):=\left\{ y \in W^{2,   p}(U ;\R^d)\,|\,
    y={\rm id} \text{ on } \Gamma_D    \right\}, 
  \end{align*}
  for almost all times    and   some given $$p>d.$$ Moreover, we impose local
  invertibility and orientation preservation.    The set
  of  
  {\it admissible deformations}    is hence   defined as
  $$ {\mathcal A} \coloneqq \left\{y \in W^{2,p}_{\Gamma_D}(U ;\R^d) \
    \Big| \  \nabla y \in \GL_+(d) \ \text{a.e. in } U\right\}.$$

\subsection{Elastic  energy} \label{sec:energy}

Let $\eps\geq 0$  be given and     $h_\eps \in
C^{\infty}(\R;[0, 1])$ for $\eps>0$ be    nondecreasing   functions such that 
\begin{equation}\label{h_eps properties}
  h_\eps( \sigma)=\begin{cases}
    0 \quad & \text{if} \ \   \sigma\leq -\eps/2,\\
    1  \quad & \text{if} \ \   \sigma\geq \eps/2,
  \end{cases}
  \quad \|   h_\eps' \|_{L^{\infty}(\R)}\leq \frac{2}{\eps}.
\end{equation} 
 Moreover, let   $ h_0$ be the  discontinuous  Heaviside-like function defined as $h_0( \sigma)=0$ if $ \sigma<
0$ and $h_0( \sigma)=1$ if $ \sigma\geq0$.  Note that 
$h_\eps \rightarrow h_0$  in  $\R\setminus \{0\}$   as $\eps\rightarrow 0$.

We  define  the elastic energy
density $W_\eps:\Rz \times \GL_+(d)\to [0,\infty) $    of the medium
  
as 
  \begin{equation}\label{W0}
    W_\eps(\sigma,F)\coloneqq (1-h_\eps(\sigma))   V^a   (F) +
    h_\eps(\sigma)    V^{ r}   (F)   + \J (F) .
  \end{equation}
 Here,  $\sigma$  is a placeholder for $\theta(x){-}t$, whose $0$-sublevel set
 $\{x\in U\mid \theta(x)<t\}$  represents  the  accreting  phase
 at time    $t>0$.    In particular, $W_\eps(\sigma, \cdot)=    V^a
     + \J$
 for $\sigma<-\eps/2$, so that $   V^a     + \J$ is the elastic energy density of the
 accreting phase. On the other hand, $W_\eps(\sigma,\cdot)=    V^r
     + \J$ for $\sigma >\eps/2$ and $    V^r     + \J$ is the elastic energy density of the
 receding phase.

 On the elastic energy densities we require 
\begin{align}
   &\label{W1}     V^a  ,\,    V^r   ,\, \J  \in C^1(\GL_+(d);[0,\infty)), \\[2mm] 
   &   \exists \cW>0\,: \quad 
     V^a(F),\, V^r(F)\geq \cW|F|^{p}-\frac{1}{\cW}
     , \nonumber\\\label{W3}
  &\qquad    V^a(F) - V^r(F) \leq \frac{1}{\cW}(1+|F|^{p})\quad \forall
     F\in \GL_+(d) , \\
   &\label{term J}   
     \exists q>\frac{pd}{p-d}\,: \quad  \J(F)\geq \frac{\cJ}{(\det F)^q}.
\end{align}
   The upper bound on $V^a - V^r $ in \eqref{W3} will be
instrumental in order to prove a control on the power associated with
the phase transformation. In particular, if the receding phase has a
higher energy density, namely, $V^r \geq V^a$, such upper bound
trivially holds.

Although    not strictly needed for the analysis we    also   require
 the frame indifference 
 \begin{equation}\label{W2}  
        V^a   (Q F)=   V^a   (F),\,    V^r   (Q F)=  
     V^r   (F)  , \, \J(QF)=\J(F)  \ \ \forall  F \in \GL_+(d), \ Q \in \SO(d).
  \end{equation}
  
As regards the second-order potential $H$ we ask for 
\begin{align} 
  &\label{H1} H \in C^1(\R^{d\times d \times
  d};[0,\infty))\ \ 
    \text{convex},\\ 
  &\label{H2} H (Q G)=H (G) \ \ \text{for all} \ \ G \in \R^{d\times d
    \times d}, \ Q \in \SO(d),\\
  &    \exists  c_H>0:  \quad c_H|G|^p\leq H(G) \leq    \frac{1}{c_H} (1+|G|)^p,\quad| {\rm D} H(G)|\leq     \frac{1}{c_H} |G|^{p-1}, \label{H3} \\
    &\qquad c_H|   G-\haz G |^p\leq ( {\rm D} H(   G )- {\rm D}
      H(   \haz G  )){\vdots}(   G-\haz G )\quad \forall   G,\, \haz G \in\R^{d \times d \times d}\label{H32} .
\end{align}
Again, the frame-indifference requirement 
\eqref{H2} is not strictly needed for the analysis. 

 By integrating over the reference configuration $U$ we define 
$\W_\eps \colon C(   \overline{U} )\times {\mathcal A}\rightarrow
[0,\infty)$ and
$\HH \colon    \mathcal{A} \rightarrow [0,\infty)$ as
\begin{equation*}
  \W_\eps(\sigma,y)\coloneqq \int_{U }W_\eps(\sigma,\nabla y)\,\dx\quad 
  \text{and}  \quad
  \HH(y)\coloneqq \int_{U }H(\nabla^2 y)\,\dx.
\end{equation*}

\subsection{Viscous  dissipation}\label{sec:dissipation} 
   For $\eps\geq 0$ given, set   the instantaneous viscous dissipation density
$R_\eps:\Rz \times \GL_+(d) \times \Rz^{d\times d}\to [0,\infty)$ as
 \begin{equation}\label{R0}
    R_\eps(\sigma,F,\dot F)\coloneqq (1-h_\eps(\sigma))R^a(F,\dot F) + h_\eps(\sigma) R^{ r}(F,\dot F) 
  \end{equation}
Here, $R^a,\, R^{ r}:\GL_+(d)  \times \R^{d\times d}\rightarrow
[0,\infty)$   are the instantaneous viscous dissipation
densities of the  accreting and of the receding phase,
respectively.    They  are    assumed to be quadratic in the rate   
$\dot C\coloneqq \dot F^\top F+
F^\top\dot F$    of
the    right   Cauchy--Green
tensor    $ C\coloneqq  F^\top F$,   namely 
\begin{equation*}
  R^a(F,\dot F)\coloneqq \frac{1}{2}\dot C {:}\DD^a(C) {:}\dot C,  \quad R^{ r}(F,\dot F)\coloneqq \frac{1}{2}\dot C {:}\DD^{ r}(C) {:}\dot C \qquad
  \forall  F \in \GL_+(d), \  \dot F\in \R^{d\times d}.
\end{equation*}
We assume that  
\begin{align}
  &\nonumber\DD^a,\, \DD^{ r}\in C(\R^{d\times d}_{\rm sym};\R^{d\times d\times
  d\times d}) \ \ \text{with}  \
    \ (\DD^i)_{jk\ell m}=(\DD^i)_{kj \ell m}=(\DD^i)_{\ell mjk} \\
  &\quad \forall j,\,k,\,\ell,\,m=1,\dots, d, \ \text{for} \
    i=a,\, r , \label{R1} \\
  &\exists c_{\DD}>0: \quad 
    c_{\DD}|\dot C|^2\leq  \dot C {:}\DD^{ i}(C) {:}\dot C \quad \forall
    C,\,\dot C\in \R^{d\times d}_{\rm sym},   \ \text{for} \ i=a,\, 
    r.
 \label{R2} 
\end{align}
Notice that  this  specific choice of $R_\eps$  ensures that
\begin{align*}
  &\partial_{\dot F}R_\eps(\sigma,F,\dot F)= 2(1{-}h_\eps(\sigma))F\DD^a(C) {:}\dot C +2h_\eps(\sigma)F\DD^{ r}(C) {:}\dot C\\
  &\quad=  2(1{-}h_\eps(\sigma))F\DD^a(F^\top F) {:} (\dot F^\top F{+} F^\top\dot F)
    + 2h_\eps(\sigma)F\DD^{ r} (F^\top F){:}(\dot F^\top F{+} F^\top\dot F),
\end{align*} 
which is of course linear in $\dot F$. 
By integrating on the reference configuration $U$ we define 
$\RR_\eps\colon C(  \overline{U} )\times    \mathcal{A}   
\times H^{1}(U
;\R^d)\rightarrow [0,\infty)$  as 
\begin{equation*}
  \RR_\eps(\sigma,y,\dot{y})\coloneqq \int_U  R_\eps(\sigma,\nabla y,\nabla \dot y)\,\dx.
\end{equation*}

\subsection{Loading and initial data}  We assume that the
body force density $f=f(\sigma,x)$ is (constant in time and) suitably
smooth with respect to $\sigma$, namely 
\begin{equation}
  \label{L1} f \in W^{1, \infty}(\R;L^{2}(U ;\R^d)),
\end{equation} 
The $\sigma$-dependence of the force density $f$ is intended to cover
the case of gravitation $f=\rho g$, where the density $\rho$ depends
on the phase, while the acceleration field $g$ is given.

We moreover assume that the initial deformation $y_0$ satisfies 
\begin{align} 
  &y_0\in \mathcal{A} \ \ \text{with} \   \int_U     V^a   (\nabla y_0 )+   V^r   (\nabla y_0 )  +\J(\nabla y_0) +H(\nabla^2 y_0)\,\dx<\infty.\label{I1}
\end{align}

\subsection{Growth}
Concerning the accretive-growth model we ask for 
\begin{align} 
  \gamma \in C^{0,1}(\Rz^d \times \GL_+(d)) \ \
  \text{with} \ \ c_\gamma\leq \gamma(\cdot) \leq C_\gamma  \
  \ \text{for some}  \ \ 0<c_\gamma\leq C_\gamma  \label{bounds
  gamma} .
\end{align}
Let moreover the initial location
 of the accreting phase be given by
\begin{equation} 
\emptyset \not =\Omega_0 \subset \subset U\ \ \text{with} \ \ \Omega_0
\ \ \text{open and} \ \ 
\label{hp dist Omega0 to partial B}
  \Omega_0+B_{C_\gamma  T}\subset \subset U. 
\end{equation}
As it will be clarified later, this last requirement guarantees that
the accreting phase does not touch the boundary $\partial U $ over the
time interval $[0,T]$,  see \eqref{eq:notouch} below.

\subsection{Main results}\label{sec:result}
 Assumptions    \eqref{h_eps properties}--\eqref{hp dist Omega0 to
   partial B}   will be 
assumed throughout the remainder of the paper. 
We are interested in solving
\eqref{eq:y_intro}--\eqref{eq:initial_intro}  in the following weak/viscosity sense.

\begin{definition}[Weak/viscosity solution]\label{def:weak sol} 
  We say that a  pair  $$(y, \theta)  \in
  \left(L^{\infty}(0,T; W^{2,p}(U;\R^d ))\cap H^1(0,T;H^1(U ;\R^d))\right) \times
      C^{0,1}(\overline U) $$
  is a
  \emph{weak/viscosity solution} to the initial-boundary-value problem
  \eqref{eq:y_intro}--\eqref{eq:initial_intro} if $y(t,\cdot) \in \mathcal A$ for all $t\in (0,T)$,  $y(0,\cdot)=y_0$, and
    \begin{align}\label{weak sol eq}
       &\int_0^T\!\! \int_{U } \big(\partial_F
         W_\eps(\theta{-}t,\nabla y){:}\nabla z{+}\partial_{\dot{F}}
         R_\eps(\theta{-}t,\nabla y, \nabla \dot{y}){:}\nabla z
        +{\rm D}H\left(\nabla^2 y\right) {\vdots} \nabla^2 z \big)\,\dx \,\dt\notag\\
       &\qquad=
         \int_0^T \!\!\int_{U } f(\theta{-}t) {\cdot}z \,\dx \,\dt \quad
\forall z\in C^{\infty}(   \overline{Q };  \R^d) \ \text{with}
         \ z=0 \ \text{on} \  \Sigma_D,
      \end{align} 
     and $\theta$ is a viscosity solution to
\begin{align}
     & \gamma\big(y(\theta(x)   \wedge T ,x),\nabla y(\theta(x)    \wedge T ,x)\big)|\nabla
     (-\theta)(x)|=1 \ \   \text{in}  \ \ U
       \setminus\overline{\Omega_0},\label{eq: eikonal system}\\ 
    &   \theta =0   \ \ \text{in}  \ \ \Omega_0. \label{eq: eikonal system2}
\end{align}
Namely, $\theta$ satisfies \eqref{eq:
  eikonal system2}, and, for all $x_0 \in U\setminus
\overline{\Omega_0}$ and any smooth function $\varphi:U \to \R$ with
$\varphi(x_0)=-\theta(x_0)$ and $\varphi\geq -\theta$ ($\varphi\leq
-\theta$, respectively) in a
neighborhood of $x_0$, it holds that $\gamma(y(\theta(x_0)    \wedge T ,x_0),
\nabla y(\theta(x_0)    \wedge T ,x_0))|\nabla
\varphi(x_0))|\leq 1 (\geq 1, \ \text{respectively})$.    Moreover,
we ask that  
\begin{equation}
  \label{eq:lip0}
  0<\frac{1}{C_\gamma} \leq |\nabla \theta| \leq \frac{1}{c_\gamma} \ \
  \text{a.e. in} \ \ U. 
\end{equation}  
\end{definition}

   Note that this weak notion of solution    in Definition
\ref{def:weak sol}   still entails the validity of
an energy equality.    Namely, we have the following.

\begin{proposition}[Energy equality]\label{prop:energy}
  Under assumptions    \eqref{h_eps properties}--\eqref{hp dist Omega0 to partial
    B},   in the
  diffused-interface case $\eps>0$ a weak/viscosity solution
  $(y,\theta)$ fulfills  for all $t\in [0,T]$ the energy equality
  \begin{align}
    &\int_U  \big(W_\eps (\theta{-}t,\nabla y )+H(\nabla^2 y)-
      f(\theta{-}t){\cdot}y\big) \,\dx  - \int_U  \left(W_\eps (\theta,\nabla y_0 )+H(\nabla^2 y_0)- f(\theta){\cdot}y_0\right) \,\dx \nonumber\\
    &\quad = -2\int_0^{t}\!\!\int_U  R_\eps\left(\theta{-}s,\nabla
      y,\nabla\dot{y}\right)\,\dx\,\ds
      -\int_{0}^{t}\!\!\int_U  \dot{f}(\theta{-}s){\cdot} y\,\dx\,\ds
      \nonumber\\
    &\qquad -\int_0^{t}\!\!\int_U \partial_\sigma W_\eps (\theta{-}s,\nabla y )\,\dx\,\ds.\label{eq:equality}
  \end{align}
 In the sharp-interface case $\eps=0$, for all $t\in [0,T]$, one has instead
 \begin{align}
   &\int_U  \big(W_0 (\theta{-}t,\nabla y )+H(\nabla^2 y)-
      f(\theta{-}t){\cdot}y\big) \,\dx  - \int_U  \left(W_0 (\theta,\nabla y_0 )+H(\nabla^2 y_0)- f(\theta){\cdot}y_0\right) \,\dx \nonumber\\
    &\quad = -2\int_0^{t}\!\!\int_U  R_0\left(\theta{-}s,\nabla
      y,\nabla\dot{y}\right)\,\dx\,\ds
      -\int_{0}^{t}\!\!\int_U  \dot{f}(\theta{-}s){\cdot} y\,\dx\,\ds
      \nonumber\\
    &\qquad -\int_0^{t}\!\!\int_{\{\theta=s\}} \frac{   V^r  (\nabla y
    )-   V^a  (\nabla y)}{|\nabla \theta|}\,{\rm d}\HH^{d-1}\,\ds.\label{eq:equality2}
  \end{align}
\end{proposition}
 
Relations \eqref{eq:equality}--\eqref{eq:equality2} express the energy
balance in the model. In particular, the difference between the actual
and the initial complementary energies (left-hand side in \eqref{eq:equality}--\eqref{eq:equality2}) equals the sum of the total
viscous dissipation, the work of external forces,
and the energy stored in the medium in connection with
the phase-transition process (the three terms in the right-hand side of
\eqref{eq:equality}--\eqref{eq:equality2}, up to    signs).   
 Proposition \ref{prop:energy} is proved in Section \ref{sec:energy_equality}.

 Our main result reads as follows. 

\begin{theorem}[Existence]  
  \label{Thm:existence sol}
   Under assumptions    \eqref{h_eps properties}--\eqref{hp dist Omega0 to partial
    B},   for all given $\eps\geq 0$   
  there exists a  weak/viscosity  solution  $(y, 
  \theta)$ of problem
  \eqref{eq:y_intro}--\eqref{eq:initial_intro}.
\end{theorem}

 A proof of Theorem \ref{Thm:existence sol} in the
diffused-interface case of  $\eps>0$ is based on an
iterative 
strategy: for given $y^k$ one finds    a   viscosity 
solution $\theta^k$ to \eqref{eq: eikonal system}--\eqref{eq: eikonal
  system2} (with $y$ replaced by $y^k$).
 Then,   given $\theta^k$ one can find  $y^{k+1}$ 
satisfying \eqref{weak sol eq} (with $\theta$ replaced by $\theta^k$).  Note that such $y^{k+1}$ may be
nonunique.
   As   the set of solutions $y$ for given $\theta$ is
generally not convex,     we   do not proceed via a fixed-point argument for 
multivalued maps  (see, e.g., \cite{Kakutani})    but rather   resort in directly proving the convergence of the
iterative procedure.    This argument is detailed in  
Section \ref{sec: Diffusive model}.  
 
Eventually, the proof of Theorem \ref{Thm:existence sol} in the sharp-interface case
$\eps=0$ will be obtained in Section \ref{Sec: eps to zero} by passing to the limit as $\eps \to 0$
along a subsequence of weak/viscosity solutions $(y_\eps,\theta_\eps)$ for
$\eps>0$. As a by-product, we have the following.

\begin{corollary}[Sharp-interface limit]\label{prop:sharp}
  Under assumptions    \eqref{h_eps properties}--\eqref{hp dist Omega0 to partial
    B},   let $(y_\eps,\theta_\eps)$  be \linebreak weak/viscosity solutions of the
  diffused-interface problem
  \eqref{eq:y_intro}--\eqref{eq:initial_intro} for $\eps>0$. Then,
  there exists a not relabeled subsequence such that
  $(y_\eps,\theta_\eps)\to (y,\theta)$    uniformly,    where $(y,\theta)$ is a
  weak/viscosity solution to  the
  sharp-interface problem    for   $\eps=0$.
\end{corollary} 

Before moving on, let us mention that the assumptions on the energy
and of the instantaneous viscous dissipation density could be generalized by not requiring the
specific forms \eqref{W0} and \eqref{R0}. In fact, one could directly assume to
be given $W_\eps =W_\eps(\sigma,F) $ and
$R_\eps=R_\eps(\sigma,F,\dot F)$ of the form
$$R_\eps(\sigma,F,\dot F) = \frac12 \dot C {:}{\mathbb D}(\sigma,C){:} \dot C$$
with $ {\mathbb D}\in C(\Rz\times \Rz_{\rm sym}^{d\times d};\Rz^{d \times d \times d \times d})
$ by suitably adapting the smoothness and coercivity
assumptions. Although the existence analysis could be carried out in this more
general situation with no difficulties, we
prefer to stick to the concrete choice of \eqref{W0} and \eqref{R0} as
it allows a more transparent distinction of the diffused- and
sharp-interface cases.

Moreover, let us point out that    admissible
deformations $y$ are presently   required to be solely {\it locally} injective, by
means of the constraint $\det \nabla y >0$. On the other hand, {\it global}
injectivity may also be enforced, in the spirit of
\cite{KromerRoubicek}, see also \cite{PalmerHealey} in the static and 
\cite{CesikGravinaKampschulte,CesikGravinaKampschulte2} in the dynamic
case. This however calls for keeping track of reaction forces due to a
possible self contact at the boundary $\Gamma_N$. From the technical
viewpoint, one would need to include an extra variable in the state
in order to model such reaction. The existence theory of Theorem
\ref{Thm:existence sol} can be extended to
cover this case, at the price of some notational intricacies. We however
prefer to avoid discussing global injectivity here, for the sake of exposition clarity.

\section{Proof of Proposition \ref{prop:energy}:
  energy equalities}\label{sec:energy_equality}

We firstly consider the diffused-interface setting    of $\eps>0$.   Let $(y,\theta)$
be a weak/viscosity solution to
\eqref{eq:y_intro}--\eqref{eq:initial_intro}.    In order to deduce
the energy equality, the Euler-Lagrange equation
\eqref{weak sol eq} should be tested by $\dot y$. This however
requires some care, as $\dot y$ is not regular enough to use it as
test function. We follow the argument of \cite{MielkeRoubicek}, based on
the validity of a chain rule for the functional $\mathcal{H}$. In
particular, we start by checking that \eqref{weak sol eq} can be
equivalently rewritten as
\begin{equation}
  \label{eq:abstract}
  \partial_2\mathcal{W}_\eps(\theta{-}t,y) +
  \partial_3\mathcal{R}_\eps(\theta{-}t,y,\dot y) +\partial
  \mathcal{H}(y) \ni \haz f  \quad \text{in} \ \
  (H^1_{\Gamma_D}(U;\R^{   d}))^*, \ \ \text{a.e. in} \ \ (0,T).
\end{equation}
Here, $(H^1_{\Gamma_D}(U;\R^{   d}))^*$ indicates  the
dual of $H^1_{\Gamma_D} (U;\R^{   d}):=\{z \in H^1 
(U;\R^{d}) \ | \ z=0 \ \text{on}
\ \Gamma_D\}$, the symbol $\partial$ denotes the (possibly partial) subdifferential
from $H^1_{\Gamma_D} (U;\R^{   d})$ to $(H^1_{\Gamma_D}(U;\R^{  
  d}))^*$ and $\haz f: (0,T)\to (H^1_{\Gamma_D}(U;\R^{  
  d}))^*$ is given by
$$\langle \haz f(t),z\rangle := \int_U f(\sigma{-}t){\cdot} z \, \dx
\quad \forall z \in H^1_{\Gamma_D} (U;\R^{   d})$$
where $\langle \cdot,\cdot \rangle$ is the duality pairing between
$(H^1_{\Gamma_D}(U;\R^{d}))^*$ and $H^1_{\Gamma_D}(U;\R^{
  d})$. Indeed, owing to the fact that $\nabla y \in L^\infty(Q)$ and
$\nabla \dot y \in L^2(Q)$ and using the regularities \eqref{W1}, \eqref{R1}, and
\eqref{L1} one can check that
\begin{align*}
  &\langle \partial_2\mathcal{W}_\eps(\theta{-}t,y) ,z\rangle = \int_U
    \partial_FW_\eps(\theta{-}t,\nabla y){:}\nabla z \quad \forall  z \in H^1_{\Gamma_D} (U;\R^{   d}),\\
  &\langle  \partial_3\mathcal{R}_\eps(\theta{-}t,y,\dot y) ,z\rangle = \int_U
    \partial_{\dot F}R_\eps(\theta{-}t,\nabla y,\nabla \dot
    y){:}\nabla z \quad \forall  z \in H^1_{\Gamma_D} (U;\R^{   d}),
\end{align*}
and that $\Sigma= \haz f - \partial_2\mathcal{W}_\eps(\theta{-}t,y) -
  \partial_3\mathcal{R}_\eps(\theta{-}t,y,\dot y) \in L^2(0,T;
  (H^1_{\Gamma_D}(U;\R^{d}))^*)$. On the other hand, using equation
  \eqref{weak sol eq}, the fact that $y\in
  L^p(0,T;W^{2,p}(\Omega;\Rz^d))$, and the  convexity \eqref{H1} of
  $H$ we get
  \begin{align*}
    \int_0^T \langle \Sigma, w-y\rangle \, \d t \stackrel{\eqref{weak
    sol eq}}{=} \int_0^T \!\!\int_U{\rm
    D}H(\nabla^2y){\vdots}\nabla^2(w-y)\, \d t \leq \int_0^T
    \big(\mathcal{H}(w)-\mathcal{H}(y)\big)\, \d t
  \end{align*}
  for all $w\in L^p(0,T;W^{2,p}(\Omega;\Rz^d))\cap L^2(0,T;
  H^1_{\Gamma_D}(U;\R^{d}))$. This in particular implies that $\Sigma
  \in \partial \mathcal{H}(y)$ a.e. in $(0,T)$, whence the abstract
  equation \eqref{eq:abstract} follows and the chain rule
  \cite[Prop.~3.6]{MielkeRoubicek} entails that $\mathcal{H}(y)\in
  W^{1,1}(0,T)$ and
  \begin{equation}
  \label{eq:chain} 
  \frac{\d}{\d t} \mathcal{H}(y) = \langle \Sigma,\dot y \rangle \quad
  \text{a.e. in} \ (0,T).
\end{equation}
Note that all terms in
  \eqref{eq:abstract} belong to $L^2(0,T;
  (H^1_{\Gamma_D}(U;\R^{d}))^*)$. One can hence test
  \eqref{eq:abstract} on $\dot y \in L^2(0,T;
  H^1_{\Gamma_D}(U;\R^{d}))$ and deduce that  
\begin{align}
  &\int_0^t \!\!\int_U \partial_F W_\eps(\theta{-}s,\nabla y){:}\nabla
  \dot y \,\dx\, \ds + \int_0^t \!\!\int_U \partial_{\dot
  F}R_\eps(\theta{-}s,\nabla y , \nabla \dot y){:} \nabla \dot y \, \dx \,
  \ds\nonumber\\
  &\quad +  \int_U H(\nabla^2y(t))\, \dx - \int_U H(\nabla^2y_0)\, \dx
      = \int_0^t \!\!\int_U
  f(\theta{-}s){\cdot} \dot y\, \dx\,\ds.\label{eq:2}
\end{align}
We readily have that
\begin{align}
&\int_U W_\eps(\theta{-}t,\nabla y)\,
  \dx -\int_U W_\eps(\theta,\nabla y_0)\,
  \dx =\int_0^t\frac{{\rm d}}{   \d s  } \int_U   \partial_F W_\eps(\theta{-}s,\nabla y)\,
  \dx\, \ds \nonumber\\
  &\quad=\int_0^t \!\!\int_U \partial_F W_\eps(\theta{-}s,\nabla y){:}\nabla
\dot y \,\dx\, \ds - \int_0^t \!\!\int_U \partial_\sigma
  W_\eps(\theta{-}s,\nabla y)\, \dx \, \ds. 
  \label{eq:derivative}
\end{align}
Moreover, it is a standard matter to check that $\partial_{\dot
  F}R_\eps(\sigma,F,\dot F){:}\dot F = 2 R_\eps(\sigma,F,\dot F)$, so
that
\begin{equation}
  \label{eq:30}
  \int_0^t \!\!\int_U \partial_{\dot
  F}R_\eps(\theta{-}s,\nabla y , \nabla \dot y){:} \nabla \dot y \, \dx \,
  \ds = 2\int_0^t \!\!\int_U  R_\eps(\theta{-}s,\nabla y , \nabla \dot y)\, \dx \,
  \ds,
\end{equation}
   whence the    energy equality \eqref{eq:equality}    in the
diffused-interface case $\eps>0$   follows    from \eqref{eq:2}.   

The proof of energy equality \eqref{eq:equality2} for the
sharp-interface case $\eps=0$ follows the same strategy, as one
can again establish \eqref{eq:2} (for $W_0$    and $R_0$ in place of
$W_\eps$ and $R_\eps$)   and
\eqref{eq:chain}. A notable difference is however in \eqref{eq:derivative}, which now
requires some extra care as $h_0$ is discontinuous. In particular, the energy equality
\eqref{eq:equality2} follows as soon as we prove that
\begin{align}
&\int_U W_0(\theta{-}t,\nabla y)\,
  \dx -\int_U W_0(\theta,\nabla y_0)\,
  \dx  \nonumber\\
  &\quad=\int_0^t \!\!\int_U \partial_F W_0(\theta{-}s,\nabla y){:}\nabla
\dot y \,\dx\, \ds - \int_0^t \!\!\int_{\{\theta=s\}} \frac{   V^r  (\nabla
    y) -    V^a  (\nabla y)}{|\nabla \theta|} \, {\rm d}
    \mathcal{H}^{d-1} \, \ds. 
  \label{eq:derivative2}
\end{align}
The remainder of this section is devoted to    check  
\eqref{eq:derivative2}.

To start with, let a nonnegative and even function $\rho \in C^\infty(\R)$ be
given with support in $[-1,1]$ and with $\int_\R\rho(s)\, \ds =1$. For $\eps>0$ we define $\rho_\eps(t) := \rho(t/\eps)/\eps$ and $\eta_\eps(t) := \int_{-1}^t\rho_\eps(s)\, \ds$ for all $t \in
\R$. As $\eta_\eps \to h_0$ in $\R\setminus\{0\}$, by letting
$$G_\eps(t)  := \int_U  \Big(    V^a  (\nabla y(t,x)) +
\eta_\eps(\theta(x){-}t) \big(    V^r  (\nabla y(t,x))  -   V^a  (\nabla
y(t,x))\big)  +\J(\nabla y(t,x))  \Big) \, \dx $$
we readily check that
\begin{equation}
  \lim_{\eps\to 0}G_\eps(t):=\int_U  W_0(\theta(x){-}t,\nabla y(t,x))\,
\dx\label{eq:G}
\end{equation}
for all $t \in [0,T]$. As $G_\eps \in H^1(0,T)$ we    can   compute its
time derivative at almost all times getting
\begin{align*}
     \frac{\d}{\dt}  G_\eps(t) &= \int_U  \Big( \partial_F    V^a  (\nabla
  y)+
\eta_\eps(\theta{-}t) \big( \partial_F   V^r  (\nabla y)  {-}\partial_F   V^a  (\nabla
  y)\big)   +\partial_F \J(\nabla y) \Big)    {:}\nabla \dot y
                     \, \dx \nonumber\\
  &\quad {}-\int_U\rho_\eps(\theta{-}t) \big(    V^r  (\nabla y)  {-}   V^a  (\nabla
y)\big)  \, \dx .
\end{align*}
By integrating in time, taking the limit $\eps \to 0$, and using
\eqref{eq:G},
one hence gets
\begin{align*}
  & \int_U  W_0(\theta{-}t,\nabla y)\,
\dx- \int_U  W_0(\theta,\nabla y_0)\,
    \dx  = \lim_{\eps\to 0} \big(G_\eps(t) - G_\eps(0) \big)=
    \lim_{\eps\to 0}\int_0^t    \frac{\d}{\d s} 
    G_\eps(s)\, \ds \nonumber\\
  &\quad = \lim_{\eps\to 0}\int_0^t\!\!\int_U \Big( \partial_F    V^a (\nabla
  y)+
\eta_\eps(\theta{-}s) \big( \partial_F   V^r  (\nabla y)  {-}\partial_F   V^a  (\nabla
    y)\big)   +\partial_F \J(\nabla y) \Big)   {:}\nabla\dot y  \, \dx \, \ds \nonumber\\
  &\qquad {}-\lim_{\eps\to 0}\int_0^t\!\! \int_U\rho_\eps(\theta{-}s) \big(    V^r  (\nabla y)  {-}   V^a  (\nabla
    y)\big)  \, \dx \, \ds \nonumber\\
  &\quad = \int_0^t \!\!\int_U \partial_F W_0(\theta{-}s,\nabla y){:}\nabla
\dot y \,\dx\, \ds -\lim_{\eps\to 0}\int_0^t\!\! \int_U\rho_\eps(\theta{-}s) \big(    V^r  (\nabla y)  {-}   V^a  (\nabla
    y)\big)  \, \dx \, \ds.
\end{align*}
In order to prove \eqref{eq:derivative2} it is hence sufficient to
check that
\begin{align}
&\lim_{\eps\to 0}\int_0^t\!\! \int_U\rho_\eps(\theta{-}s) \big(    V^r  (\nabla y)  {-}   V^a  (\nabla
  y)\big)  \, \dx \, \ds  =\int_0^t \!\!\int_{\{\theta=s\}} \frac{   V^r  (\nabla
    y) {-}    V^a  (\nabla y)}{|\nabla \theta|} \, {\rm d}
    \mathcal{H}^{d-1} \, \ds.
  \label{eq:tocheck}
\end{align}
By introducing the short-hand notation $g =    V^r  (\nabla
  y) -    V^a  (\nabla y)$ and by
using the coarea formula \cite[Sec. 3.2.11]{Federer} (recall that
$\theta$ is Lipschitz continuous and $|\nabla \theta| \geq
1/C_\gamma   >0  $
almost everywhere,    see \eqref{eq:lip0})   we can compute
\begin{align}
  &\int_0^t\!\! \int_U\rho_\eps(\theta{-}s) \big(    V^r  (\nabla y)  {-}   V^a  (\nabla
  y)\big)  \, \dx \, \ds = \int_0^t\!\!\int_\R \int_{\{\theta=r\}}
    \rho_\eps(\theta   (x)  {-}s)\frac{g(s,x)}{|\nabla \theta(x)|} {\rm d}\mathcal{H}^{d-1}(x)\, 
    {\rm d}r \, \ds\nonumber\\
  &\quad = \int_0^t\!\!\int_\R \int_{\{\theta=r\}}\rho_\eps(r{-}s)
    \frac{g(r,x)}{|\nabla \theta(x)|} {\rm d}\mathcal{H}^{d-1}(x)\, 
    {\rm d}r \, \ds \nonumber\\
  &\qquad + \int_0^t\!\!\int_\R \int_{\{\theta=r\}}\rho_\eps(r{-}s)
    \frac{g(s,x){-}g(r,x)}{|\nabla \theta(x)|} {\rm d}\mathcal{H}^{d-1}(x)\, \
    {\rm d}r \, \ds. \label{eq:rhs}
\end{align}
The coarea formula and the bound $|\nabla \theta|\leq 1/c_\gamma$   
(see again \eqref{eq:lip0})  
ensure that $r\in \R \mapsto    m  (r):=
\mathcal{H}^{d-1}(\{\theta = r\})$ is integrable. Indeed,
$$\|    m   \|_{L^1(\R)}=\int_\R \mathcal{H}^{d-1}(\{\theta = r\}) \, {\rm d} r= \int_U
|\nabla \theta |\,\dx <\infty.$$
As $g/|\nabla \theta|$ is bounded, setting 
$$r\in \Rz \mapsto\ell(r) := \int_{\{\theta=r\}}
\frac{g(r,x)}{|\nabla \theta(x)|} {\rm d}\mathcal{H}^{d-1}(x)$$
one has that $\ell \in L^1(\R)$, as well,    since
\begin{align*}
  \| \ell \|_{L^1(\Rz)} = \int_\Rz\!\! \int_{\theta = r}
  \frac{|g(r,x)|}{|\nabla \theta(x)|}\,\d \mathcal{H}^{d-1}(x)\, \d r
  \leq \sup\frac{|g|}{|\nabla \theta|} \| m \|_{L^1(\Rz)}<\infty.
\end{align*}
Moreover, we have that  
\begin{align*}& \int_\R \int_{\{\theta=r\}}\rho_\eps(r{-}s)
    \frac{g(r,x)}{|\nabla \theta(x)|} {\rm d}\mathcal{H}^{d-1}(x)\, 
                {\rm d}r    \\
  &\quad = \int_\R \rho_\eps(s{-}r)\left(\int_{\{\theta=r\}}
    \frac{g(r,x)}{|\nabla \theta(x)|} {\rm d}\mathcal{H}^{d-1}(x)\right) \,
    {\rm d}r     =  (\rho_\eps \ast \ell)(s)
    \end{align*}
 where we used that $\rho_\eps$ is even and where the symbol $\ast$ stands
 for the    usual   convolution in $\R$. As
 $\rho_\eps\ast \ell \to \ell$ strongly in $L^1(\R)$ for  $\eps \to 0$, 
one can pass to the limit  in the  
first term on the right-hand side of \eqref{eq:rhs} and get
\begin{align}
  &\lim_{\eps \to 0}\int_0^t\!\!\int_\R \int_{\{\theta=r\}}\rho_\eps(r{-}s)
    \frac{g(r,x)}{|\nabla \theta(x)|} {\rm d}\mathcal{H}^{d-1}(x)\, 
    {\rm d}r \, \ds  = \lim_{\eps \to 0} \int_0^t (\rho_\eps\ast \ell)\,\ds
  = \int_0^t \ell \, \ds \nonumber\\
  &\quad = \int_0^t\!\!\int_{\{\theta =s\}}
  \frac{g(s,x)}{|\nabla \theta(x)|} \,{\rm d}\mathcal{H}^{d-1}(x)\,   
    {\rm d}s = \int_0^t\!\!\int_{\{\theta =s\}}
  \frac{   V^r  (\nabla y){-}   V^a  (\nabla y)}{|\nabla \theta|}\, {\rm d}\mathcal{H}^{d-1}\,   
    {\rm d}s.
  \label{eq:40}  
\end{align}
As regards the second term in the right-hand side of \eqref{eq:rhs}, 
notice that $\rho_\eps(r{-}s)\neq 0$ only if $|r-s|\leq 2\eps$. Hence, using 
 the H\"older regularity of  $g$ and the boundedness of
$1/|\nabla \theta|$ and $|\nabla \theta|$, we conclude that
\begin{align}
  \label{eq:10}
  &\lim_{\eps \to 0}\left|\int_0^t\!\!\int_\R \int_{\{\theta=r\}}\rho_\eps(r{-}s)
    \frac{g(s,x){-}g(r,x)}{|\nabla \theta(x)|} {\rm d}\mathcal{H}^{d-1}(x)\, 
  {\rm d}r \, \ds\right| \nonumber\\
  &\quad \leq \lim_{\eps \to 0}    c  \,\eps^\alpha  \int_0^t\!\!\int_\R  \rho_\eps(r{-}s)  \,\mathcal{H}^{d-1}(\{\theta=r\})\, 
    {\rm d}r \, \ds \nonumber\\
  &\quad = \lim_{\eps \to 0}    c  \, \eps^\alpha \|
    \rho_\eps\ast    m   \|_{L^1(\Rz)}\leq\lim_{\eps \to 0}    c  \, \eps^\alpha \|
        m \|_{L^1(\Rz)} =0
\end{align}
for some $\alpha\in (0,1)$.
Relations \eqref{eq:40}--\eqref{eq:10} imply that the limit
\eqref{eq:tocheck} holds true. This in turn proves
\eqref{eq:derivative2} and the energy equality \eqref{eq:equality2}
follows.

\section{Proof of Theorem \ref{Thm:existence
  sol}: diffused-interface case}\label{sec: Diffusive model}

 Let $\eps>0$ be fixed. We prove existence of a weak/viscous
solution $(y,\theta)$ by an iterative construction.    We
  start     by proving that for all given $\theta\in C(   \overline{U}  )$ there
exists  an admissible  deformation 
$y$ 
satisfying \eqref{weak sol eq}.

\begin{proposition}[Existence of $y$ given $\theta$]\label{prop: existence y}
  Set $\eps> 0$ and let
  $\theta\in C(   \overline{U}  )$ be fixed    with $\Omega(T)
  \subset\subset U$.     Under assumptions   
  \eqref{h_eps properties}--\eqref{I1}    there
  exists $y\in L^{\infty}(0,T;W^{2,p}(U ;\R^d))\cap
  H^1(0,T;H^1(U ;\R^d))$  with $ y(t,\cdot) \in \mathcal{A}$ for
  every $t\in (0,T)$ satisfying \eqref{weak sol eq}.    More
  precisely, there exists a positive constant $c$ depending on data but independent of
  $\eps$ and $\theta$ such that
  \begin{equation}
    \label{eq:exbound}
    \| y\|_{L^{\infty}(0,T;W^{2,p}(U ;\R^d))\cap
  H^1(0,T;H^1(U ;\R^d))}\leq c.  
  \end{equation} 
\end{proposition}

\begin{proof}    
The assertion follows by  adapting the arguments in
\cite{BadalFriedrichKruzik} or \cite{KromerRoubicek}.   
   We proceed by time discretization.    
Let  the time step $\tau\coloneqq
T/N_\tau$  with  $N_\tau \in \N$ be given   and let $t_i\coloneqq i \tau $, for
$i=0,\dots, N_\tau$ be the corresponding uniform partition of the time
interval $[0,T]$.    Within this  proof, the generic constant $c$ is always independent of
  the given $\theta$, as well.  

 For $i=1,\dots, N_\tau$, we  define 
$y^i_\tau\in  \mathcal{A}$     via  
  \begin{align*}
    y^i_\tau\in \argmin_{ y\in\mathcal A}\left\{\W_\eps(\theta{-}t_i,y)+\HH(y)+\tau \RR_\eps\left(\theta{-}t_i,y^{i-1}_\tau,\frac{y{-}y^{i-1}_\tau}{\tau}\right) -\int_U  f(\theta{-}t_i){\cdot}y\,\dx \right\}.
  \end{align*}
Under the growth conditions \eqref{W3}--\eqref{term J},   \eqref{H3},  and
 \eqref{R2},    and   the regularity    and convexity   assumptions \eqref{W1}, \eqref{H1},
\eqref{R1},  and  \eqref{L1}, the existence of $y^i_\tau$ for
every $i=1,\dots, N_\tau$ follows by the Direct Method of the calculus
of variations.  Moreover, every minimizer $y_\tau^i$ satisfies the time-discrete
Euler--Lagrange equation

\begin{align}\label{Euler-Lagrange discrete}
  &\int_{U } \left(\partial_F W_\eps(\theta{-}t_i,\nabla
    y^i_\tau){+}\partial_{\dot{F}} R_\eps\left(\theta{-}t_i,\nabla
    y^{i-1}_\tau, \frac{\nabla y^i_\tau{-}\nabla
    y^{i-1}_\tau}{\tau}\right)\right){:}\nabla    z^i    \,\dx \notag\\
  &\quad+\int_U  { \rm D}H\left(\nabla^2 y^i_\tau\right) {\vdots}
    \nabla^2    z^i    \,\dx =
  \int_U  f(\theta{-}t_i){\cdot}z\,\dx  
 \end{align}
for every    $z^i \in \mathcal{A} $.    

Let us introduce the following notation for the time interpolants on
the partition: Given a vector $(u_{0},...,u_{N_\tau})$, 
we define its backward-constant interpolant $\overline{u}_{\tau}$, its forward-constant interpolant $\underline{u}_{\tau}$, and its piecewise-affine interpolant $\haz{u}_{\tau}$ on the partition $(t_i)_{i=0}^{N_\tau}$ as
  \begin{align*}
    &\overline{u}_\tau(0):=u_0, \quad \quad \overline{u}_\tau(t):=u_{i}  &\text{ if } t\in(t_{i-1},t_i] \quad\text{ for } i=1,\dots, N_\tau,\\
    &\underline{u}_{\tau}(T):=u_{N_\tau}, \quad \; \underline{u}_\tau(t):=u_{i-1} &\text{ if } t\in[t_{i-1},t_i) \quad\text{ for } i=1,\dots, N_\tau,\\
    &\haz{u}_\tau(0):=u_0, \quad \quad \haz{u}_\tau(t):=\frac{u_i-u_{i-1}}{t_i-t_{i-1}}(t-t_{i-1})+u_{i-1}
    &\text{ if } t\in(t_{i-1},t_i] \quad\text{ for } i=1,\dots, N_\tau.
  \end{align*}
   Owing to this notation, we can take the sum in
  \eqref{Euler-Lagrange discrete} for $i=1,\dots,N_\tau$ and equivalently
  rewrite    the discrete Euler--Lagrange equation   in the compact form 
\begin{align}
  &\int_0^T\!\!\int_{U } \left(\partial_F W_\eps(\theta{-}\overline{t}_\tau,\nabla
    \overline{y}_\tau){+}\partial_{\dot{F}} R_\eps\left(\theta{-}\overline{t}_\tau,\nabla
    \underline{y}_\tau, \nabla \dot{\haz{y}}_\tau
    \right)\right){:}\nabla    \overline{z}_\tau    \,\dx\, \dt \notag\\
  &\quad+\int_0^T\!\!\int_U  { \rm D}H\left(\nabla^2 \overline{y}_\tau\right)
    {\vdots} \nabla^2    \overline{z}_\tau  \,\dx \, \dt=
  \int_0^T\!\!\int_U  f(\theta{-}\overline{t}_\tau){\cdot}  
    \overline{z}_\tau  \,\dx \, \dt \label{Euler-Lagrange compact} 
 \end{align}
     where $\overline{z}_\tau$ is the    backward-constant    interpolant
  of $(z^i)_{i=1}^{N_\tau}$.  

  From the minimality of $y^i_\tau$ we get that 
\begin{align}
  &\int_U  W_\eps (\theta{-}t_i,\nabla y^i_\tau )\, \dx
    +\int_UH(\nabla^2 y^i_\tau)\, \dx -\int_U  f(\theta{-}t_i){\cdot}y^{i}_\tau\,\dx \nonumber\\
  &\qquad +\tau\int_U R_\eps\left(\theta{-}t_i,\nabla
    y^{i-1}_\tau,\frac{\nabla y^i_\tau{-}\nabla y^{i-1}_\tau}{\tau}\right)\,\dx   \nonumber \\
  &\quad \leq \int_U  W_\eps (\theta{-}   t_{i-1} ,\nabla y^{i-1}_\tau )\,
    \dx+\int_UH(\nabla^2 y^{i-1}_\tau)\,\dx-\int_U
    f(\theta{-}   t_{i-1} ){\cdot}y^{i-1}_\tau\,\dx\nonumber\\
  &\qquad   -\int_U
    \big(f(\theta{-} t_{i})- f(\theta{-}
    t_{i-1})\big){\cdot}y^{i-1}_\tau\,\dx\nonumber\\
  &\qquad   - \int_U
    \big(W_\eps (\theta{-} t_{i-1},\nabla y^{i-1}_\tau ) -W_\eps (\theta{-} t_{i},\nabla y^{i-1}_\tau 
    )\big) \,
    \dx. \label{eq:perstime}
\end{align}
  
Summing over $i=1,\dots,n\leq N_\tau$ inequality
\eqref{eq:perstime} we get 
\begin{align}
  &\int_U  W_\eps (\theta{-}t_n,\nabla y^n_\tau   )\, \dx+\int_U  H(\nabla^2
    y^n_\tau)\,\dx-\int_U  f(\theta{-}t_n){\cdot}y^n_\tau\,\dx \notag\\
  &\qquad+\sum_{i=1}^{n}\tau\int_U R_\eps\left(\theta{-}t_i,\nabla
    y^{i-1}_\tau,\frac{\nabla y^i_\tau{-}\nabla
    y^{i-1}_\tau}{\tau}\right)\,\dx   \nonumber \\
  &\quad \leq \int_U  W_\eps (\theta ,\nabla y_0 )+H(\nabla^2
    y_0)\,\dx-\!\int_U \! f(\theta){\cdot} y_0\,\dx\notag\\
  &\qquad-\sum_{i=1}^n\int_U
    (f(\theta{-}t_i){-}f(\theta{-}t_{i-1})){\cdot}y^{i-1}_\tau\,\dx\nonumber\\
  &\qquad   -\sum_{i=1}^n\int_U
    \big(W_\eps (\theta{-} t_{i-1},\nabla y^{i-1}_\tau ) -W_\eps (\theta{-} t_{i},\nabla y^{i-1}_\tau 
    )\big) \,
    \dx.\label{energy ineq 10}
\end{align}
By the growth conditions \eqref{W3}--\eqref{term J} \eqref{H3},
\eqref{R2}, and \eqref{L1}, we hence have that 
\begin{align}
  &   \cW \|\nabla y^n_\tau\|_{L^{p}(U; \Rz^{d\times d})}^p
    +   \cJ \left\|\frac{1}{\det \nabla
     y^n_\tau}\right\|^q_{L^q(U
     )}+   c_H   \|\nabla^2y^n_\tau\|^p_{L^p(U ;
    \Rz^{d\times  d\times d} )}    \notag\\
  &\qquad   + c_{\mathbb{D}} \sum_{i=1}^n\tau\int_U \left|\frac{(\nabla
    y^i_\tau  - \nabla y^{i-1}_\tau)^\top}{\tau} \nabla
    y^{i-1}_\tau + (\nabla y^{i-1}_\tau)^\top \frac{\nabla
    y^i_\tau - \nabla y^{i-1}_\tau}{\tau}  \right|^2\, \dx - \frac{|U|}{\cW}   \notag\\
  &\quad\leq 
    \int_U  W_\eps (\theta{-}t_n,\nabla y^n_\tau   )\, \dx+\int_U  H(\nabla^2
    y^n_\tau)\,\dx+\sum_{i=1}^{n}\tau\int_U R_\eps\left(\theta{-}t_i,\nabla
    y^{i-1}_\tau,\frac{\nabla y^i_\tau{-}\nabla
    y^{i-1}_\tau}{\tau}\right)\,\dx 
  \notag\\
&\quad \hspace{-1.5mm}\stackrel{\eqref{energy ineq 10}}{\leq} \int_U  W_\eps (\theta ,\nabla y_0 )+H(\nabla^2
    y_0)\,\dx-\!\int_U \! f(\theta){\cdot} y_0\,\dx+\int_U  f(\theta{-}t_n){\cdot}y^n_\tau\,\dx \notag\\
  &\qquad   +    \sum_{i=1}^n \int_{t_{i-1}}^{t_i}\int_U
    \dot f(\theta{-}s){\cdot}y^{i-1}_\tau\,\dx\, \d s \nonumber\\
  &\qquad {}-  \sum_{i=1}^n\int_U
    \big(W_\eps (\theta{-} t_{i-1},\nabla y^{i-1}_\tau ) -W_\eps (\theta{-} t_{i},\nabla y^{i-1}_\tau 
    )\big) \,
    \dx . \label{eq:usamidopo} 
\end{align}  
In order to control the right-hand side above, we remark that 
\begin{align*}
  &{}-  \sum_{i=1}^n\int_U
    \big(W_\eps (\theta{-} t_{i-1},\nabla y^{i-1}_\tau ) -W_\eps (\theta{-} t_{i},\nabla y^{i-1}_\tau 
    )\big) \, \dx \\
  &\quad \hspace{-1mm} \stackrel{\eqref{W0}}{=} \sum_{i=1}^n\int_U \big(h_{   \eps}(\theta{-} t_{i-1})-h_{   \eps}(\theta{-}
    t_{i})\big)\big(V^a(\nabla y^{i-1}_\tau) - V^r(\nabla
    y^{i-1}_\tau) \big)\, \d x\\
  &\quad\hspace{-1mm} \stackrel{\eqref{W3}}{\leq}  \frac{1}{\cW} \sum_{i=1}^n\int_U \big(h_{   \eps}(\theta{-} t_{i-1})-h_{   \eps}(\theta{-}
    t_{i})\big)\big(1+|\nabla
    y^{i-1}_\tau|^p \big)\, \d x
\end{align*}  
where we have also used that 
$h_{   \eps}(\theta{-} t_{i-1})-h_{   \eps}(\theta{-}
t_{i} )\geq 0$. As $h_{   \eps}(\theta{-} t_{i-1})-h_{   \eps}(\theta{-}
t_{i} )= 0$ on the complement of $$E_i:=\{x\in U \:: \: \theta(x) \in
[t_{i-1}-\eps/2,t_i+\eps/2]\},$$
by using 
    $\|h'_\eps\|_{L^\infty(\Rz)}\leq 2/\eps$    (recall \eqref{h_eps
    properties})    and the embedding $L^\infty(U,\Rz^{d\times d})
    \subset W^{1,p}(U,\Rz^{d\times d})$  we get  
\begin{align*}
  &{}-  \sum_{i=1}^n\int_U
    \big(W_\eps (\theta{-} t_{i-1},\nabla y^{i-1}_\tau ) -W_\eps (\theta{-} t_{i},\nabla y^{i-1}_\tau 
    )\big) \, \dx \\
  &\quad \leq \frac{c\tau}{\eps}\sum_{i=1}^n|E_i| \left(1+\|\nabla
    y^{i-1}_\tau\|^p_{L^\infty(U;\Rz^{d\times d})}\right) \\
    &\quad \leq \frac{c\tau}{\eps}\sum_{i=1}^n|E_i| \left(1+\|\nabla
   y^{i-1}_\tau\|^p_{L^p(U;\Rz^{d\times d})} + \|\nabla^2
   y^{i-1}_\tau\|^p_{L^p(U;\Rz^{d\times d \times d})}\right).
\end{align*}
Together with \eqref{L1}--\eqref{I1}, this allows to deduce from
inequality \eqref{eq:usamidopo} that
\begin{align}
&\cW \|\nabla y^n_\tau\|^{   p}_{L^{   p}(U; \Rz^{d\times d}
     )}+   \cJ \left\|\frac{1}{\det \nabla
     y^n_\tau}\right\|^q_{L^q(U
     )}+   c_H   \|\nabla^2y^n_\tau\|^p_{L^p(U ;
    \Rz^{d\times  d\times d} )}    \notag\\
  &\qquad   + c_{\mathbb{D}} \sum_{i=1}^n\tau\int_U \left|\frac{(\nabla
    y^i_\tau  - \nabla y^{i-1}_\tau)^\top}{\tau} \nabla
    y^{i-1}_\tau (\nabla y^{i-1}_\tau)^\top \frac{\nabla
    y^i_\tau - \nabla y^{i-1}_\tau}{\tau}  \right|^2\, \dx\nonumber\\
  &\quad \leq c + c\| y^n_\tau\|_{L^2(U;\Rz^d)} + c\sum_{i=1}^n\tau \|
    y^{i-1}_\tau\|_{L^2(U;\Rz^d)} \nonumber\\
  &\qquad +
    \frac{c\tau}{\eps}\sum_{i=1}^{n} |E_i|\left(1+ \|\nabla y^{i-1}_\tau\|^p_{L^p(U;\Rz^{d\times
    d})}+ \|\nabla^2 y^{i-1}_\tau\|^p_{L^p(U;\Rz^{d\times
    d\times d})}\right). \label{eq:toberefined2}
\end{align}
For $\tau<\eps$ one has that $\cup_{i=1}^{N_\tau} E_i$ covers
$\Omega(T)$ multiple times. In particular, we have that
\begin{equation}
  \label{eq:E}
  \sum_{i=1}^{N_\tau}|E_i| \leq \left(\frac{\eps+\tau}{\tau} +1\right)|\Omega(T)|.
\end{equation}
Hence, by the Poincar\'e inequality and the Discrete Gronwall Lemma
\cite[(C.2.6), p. 534]{Kruzik Roubicek} we find the bound
\begin{align}
  &\max_n \| y^n_\tau\|_{W^{2,p}(U,\Rz^d)}^{   p} +\sum_{i=1}^{N_\tau}\tau
\left\| \frac{(\nabla
    y^i_\tau  {-} \nabla y^{i-1}_\tau)^\top}{\tau} \nabla
    y^{i-1}_\tau + (\nabla y^{i-1}_\tau)^\top \frac{\nabla
    y^i_\tau {-} \nabla y^{i-1}_\tau}{\tau}  \right\|^2_{L^2(U;\Rz^{d\times d})} \nonumber\\
  &\quad \leq c \, {\rm exp} \left( \frac{c\tau}{\eps}
    \sum_{i=i}^{N_\tau}|E_i| \right)\stackrel{\eqref{eq:E}}{\leq} c \, {\rm exp} \left( \frac{c\tau}{\eps}
    \left(\frac{\eps + \tau}{\tau}+1\right)|\Omega(T)|\right) \leq c  \, {\rm exp}(c\tau/\eps), \label{eq:prebound0}
\end{align}
where we also used the fact that $\Omega(T)    \subset   \subset U$. 

By  
 the Sobolev embedding of $W^{2,p}(U ;\R^d)$ into $C^{1-d/p}(U
;\R^d)$ and the classical result of \cite[Thm. 3.1]{HealeyKromer} we get
\begin{equation}\label{determinant interpolant>0}
   \det\nabla \overline{y}_\tau\geq     c_\eps   >0 \ \ \text{in} \ \
 [0,T]\times \overline{U}
\end{equation}
   where the constant $c_\eps$ depends on the bound in
\eqref{eq:prebound0}.  

By the Poincar\'e inequality and the generalization of Korn's first
inequality by \cite{Neff} and \cite[Thm. 2.2]{Pompe},    also using
\eqref{determinant interpolant>0}   we have that 
\begin{equation*}
  \|\nabla \dot{\haz y}_\tau\|^{   2}_{   L^2(0,T;   L^2(
    Q;\R^{d\times d} ))}\leq    c_\eps' \int_0^T\ \| \nabla\dot{\haz{y}}_\tau^\top
  \nabla\underline{y}_\tau + \nabla
  \underline{y}_\tau^\top\nabla\dot{\haz{y}}_\tau  
  \|^2_{L^2(U;\Rz^{d\times d})} \, \d s \stackrel{\eqref{eq:prebound0}}{\leq}c_\eps' c  \, {\rm exp}(c\tau/\eps) 
\end{equation*} 
   where the constant $c_\eps'$ depends on the bound
\eqref{eq:prebound0} and on the constant $c_\eps$ in
\eqref{determinant interpolant>0}.   Again by the Poincar\'e inequality, this time applied to
$\dot y$, we get that 
\begin{equation}\label{bound H^1(H^1)}
  \| {\haz y}_\tau\|_{H^1(0,T;H^1(U; \R^d ))}\leq    c_\eps' c  \, {\rm exp}(c\tau/\eps).  
\end{equation} 
By using  these  estimates,    as
$\tau \to 0$,   up to not relabeled subsequences we get
\begin{align}
  &\overline{y}_\tau,\,\underline{y}_\tau \stackrel{*}{\rightharpoonup}
  y  \quad \text{weakly-$*$ in} \ \  L^{\infty}(0,T;W^{2,p}(U ;\R^d)),\label{convergence 1}\\
  &\nabla \dot{\haz y}_\tau \rightharpoonup \nabla \dot{y}
          \quad\text{weakly in} \ \  L^{2}( Q;\R^d),\label{convergence 2}\\
  &\nabla \haz y_\tau\rightarrow \nabla y \quad\text{strongly in}
    \ \  C^{0,\alpha}( \overline{Q} ;\R^d)\label{convergence 3} 
\end{align} 
for some $\alpha\in(0,1)$. 
   In particular,   from the convergences above we also get  $\det \nabla
\overline{y}_\tau \to \det \nabla y$    uniformly.   In combination
with the lower bound \eqref{determinant interpolant>0}, this implies
that $\nabla    y  
\in \GL_+(d)$ everywhere, hence $y$  is admissible, namely, $ y(t,\cdot) \in \mathcal{A}$ for every $t\in (0,T)$. 

We now pass to the limit in the time-discrete Euler--Lagrange equation
\eqref{Euler-Lagrange compact}.    Let $z\in
C^\infty(\overline{Q};\Rz^d)$ with $z=0$ on $\Sigma_D$ be given and
let $(z^i_\tau)_{i=1}^{N_\tau}\in \mathcal{A}$ be such that $\overline{z}_\tau \to z$
strongly in $L^\infty(0,T;W^{2,p}(U;\Rz^d)) $.  
By \eqref{L1} we have 
\begin{equation}
  \int_0^T\!\!\int_U  f(\theta{-}\overline{t}_\tau){\cdot}  
  \overline{z}_\tau  \,\dx
  \,\dt\rightarrow \int_0^T\!\!\int_U  f(\theta{-}t){\cdot}z\,\dx\,\dt.\label{eq:uno}
\end{equation}

   As $h_\eps(\theta(x){-}\overline{t}_\tau(t))\to
h_\eps(\theta(x){-}t)$ for almost every $(t,x)\in Q$, the  
dissipation    term    goes to the limit as follows
\begin{align}
  &\int_0^T\!\!\int_U \partial_{\dot{F}} R_\eps\left(\theta{-}\overline{t}_\tau,\nabla \underline{y}_\tau,\nabla\dot{\haz y }_\tau\right){:}\nabla   
  \overline{z}_\tau    \,\dx\,\dt\nonumber\\
  &\quad=2\int_0^T\!\!\int_U  (1{-}h_\eps(\theta{-}\overline{t}_\tau))\nabla \underline{y}_\tau  \left(\DD^{ a}(\nabla \underline{y}_\tau^\top \nabla \underline{y}_\tau)(\nabla\dot{\haz y }_\tau^\top \nabla \underline{y}_\tau{+}\nabla \underline{y}_\tau^\top \nabla\dot{\haz y }_\tau)\right){:}\nabla   
  \overline{z}_\tau  \,\dx\,\dt\nonumber\\
  &\quad\quad +2\int_0^T\!\!\int_U  h_\eps(\theta{-}\overline{t}_\tau)\nabla \underline{y}_\tau  \left(\DD^r(\nabla \underline{y}_\tau^\top \nabla \underline{y}_\tau)(\nabla\dot{\haz y }_\tau^\top \nabla \underline{y}_\tau{+}\nabla \underline{y}_\tau^\top \nabla\dot{\haz y }_\tau)\right){:}\nabla   
  \overline{z}_\tau  \,\dx\,\dt\nonumber\\ 
  &\quad\rightarrow 2\int_0^T\!\!\int_U  (1{-}h_\eps(\theta{-}t))\nabla y  \left(\DD^{ a}(\nabla y^\top \nabla y)(\nabla\dot{y}^\top \nabla y{+}\nabla y^\top \nabla\dot{y})\right){:}\nabla z\,\dx\,\dt\nonumber\\
  &\quad\quad+ 2\int_0^T\!\!\int_U h_\eps(\theta{-}t)\nabla y  \left(\DD^r(\nabla y^\top \nabla y)(\nabla\dot{y}^\top \nabla y{+}\nabla y^\top \nabla\dot{y})\right){:}\nabla z\,\dx\,\dt\\
  &\quad=\int_0^T\!\!\int_U \partial_{\dot{F}} R_\eps\left(\theta{-}t,\nabla y,\nabla\dot{y}\right){:}\nabla z  \,\dx\,\dt\label{eq:due}
\end{align}
   where we used \eqref{R1}  and   convergences \eqref{convergence 1}--\eqref{convergence 3}. Moreover, we also have 
\begin{equation}
  \int_0^T\int_{U }\partial_F W_\eps(\theta{-}\overline{t}_\tau,\nabla
  \overline{y}_\tau)   {:} \nabla \overline{z}_\tau   \,\dx\,\dt\rightarrow \int_0^T\int_{U }\partial_F W_\eps(\theta{-}t,\nabla y)    {:} \nabla {z}   \,\dx\,\dt\label{eq:tre}
\end{equation}
by \eqref{W1}    and convergences \eqref{convergence 1} and
\eqref{convergence 3}.  

For the convergence of the second-gradient term  we 
reproduce in this setting the argument from 
\cite{KromerRoubicek}.    Given the limit $y$, let
$(w_{\tau}^i)_{i=1}^{N_\tau}\in \mathcal{A}$ be such that 
$\overline{w}_{\tau} \to y$ strongly in
$L^\infty(0,T;W^{2,p}(U;\Rz^d))$.    We consider   
the   test functions $   \overline{z}_{\tau}
  \coloneqq
   \overline{w}_{\tau} -\overline{y}_{\tau}$  
in the time-discrete Euler--Lagrange equation \eqref{Euler-Lagrange compact}.
   Convergences \eqref{convergence 1}--\eqref{convergence 2} entail
that    $   \overline{z}_{\tau} \rightarrow 0$ strongly
in  $ L^\infty (0,T;H^{1}(U;\R^d))$ and $  
\overline{z}_{\tau} \weak 0$    weakly-$*$   in $ L^{  
  \infty } (0,T;W^{2,p} (U;\R^d))$.  
   Let us now compute
\begin{align}
 & \int_0^T\!\!\int_U({ \rm D}H(\nabla^2 y)-{ \rm D} H(\nabla^2
  \overline{y}_\tau)){\vdots}(\nabla^2 y -\nabla^2 
  \overline{y}_\tau)\,\dx\,\dt\nonumber\\
  &\quad = \int_0^T\!\!\int_U({ \rm D}H(\nabla^2
  y)-{ \rm D} H(\nabla^2 \overline{y}_\tau)){\vdots}(\nabla^2 
    y - \nabla^2 \overline{w}_{\tau})\,\dx\,\dt\nonumber\\
  &\qquad +\int_0^T\!\!\int_U({ \rm D}H(\nabla^2 y)-{ \rm D} H(\nabla^2 \overline{y}_\tau)){\vdots}\nabla^2  \overline{z}_{\tau}\,\dx\,\dt.\label{eq:sommadidue}
\end{align}
As $\nabla^2 \overline{w}_{\tau} \to \nabla^2 y$ strongly in
$L^{p}(   Q  ;\Rz^{d\times d \times d})$ and ${\rm D} H (\nabla^2
\overline{y}_\tau)$ is bounded in $L^{p'}(   Q  ;\Rz^{d\times d \times
  d})$ by \eqref{H3}, the first integral in the right-hand side above converges to $0$
as $\tau \to 0$. Hence, passing to the $\limsup$ in \eqref{eq:sommadidue},  
by the Euler--Lagrange equation \eqref{Euler-Lagrange compact}   
and convergences \eqref{convergence 3}--\eqref{eq:tre}   
we find that
\begin{align}
  &   \limsup_{\tau\rightarrow 0} \int_0^T\!\!\int_U({ \rm D}H(\nabla^2
    y)-{ \rm D} H(\nabla^2
    \overline{y}_\tau)){\vdots}(\nabla^2y - \nabla^2 \overline{y}_\tau)
     \,\dx\,\dt\nonumber\\
  &\quad    =   \limsup_{\tau\rightarrow 0} \int_0^T\!\!\int_U({ \rm D}H(\nabla^2 y)-{ \rm D} H(\nabla^2 \overline{y}_\tau)){\vdots}\nabla^2    \overline{z}_{\tau} \,\dx\,\dt\nonumber\\
  &\quad=\limsup_{\tau\rightarrow 0}\Bigg(\int_0^T\!\!\int_U{ \rm D}
    H(\nabla^2 y){\vdots}\nabla^2
      \overline{z}_{\tau}\, \dx \, \dt - \int_0^T\!\!\int_U  f(\theta{-}\overline{t}_\tau){\cdot}   \overline{z}_{\tau} \,\dx\,\dt
  \nonumber\\
  &\qquad+\int_0^T\!\!\int_U\left(\partial_F
    W_\eps(\theta{-}\overline{t}_\tau, \nabla
    \overline{y}_\tau)+\partial_{\dot
    F}R_\eps\left(\theta{-}\overline{t}_\tau, \nabla
    \underline{y}_\tau, \nabla
    \dot{\haz{y}}_\tau\right)\right){:}\nabla    \overline{z}_{\tau} \,\dx \,\dt \label{eq:prima}
    \Bigg)=0
\end{align}  
   The coercivity \eqref{H3} then 
implies    that    
\begin{equation*}
  \nabla^2 \overline{y}_\tau\rightarrow  \nabla^2 y \quad
  \text{strongly in} \ \  L^p(Q;   \R^{d\times d\times d} )
\end{equation*}  
and thus
\begin{equation*}
  { \rm D} H(\nabla^2 \overline{y}_\tau)\rightarrow { \rm D}
  H(\nabla^2 {y}) \quad \text{ strongly in } L^{p'}(Q;  
  \R^{d\times d \times d} ).
\end{equation*} 
Passing    to    the limit    as $\tau \to 0 $    in \eqref{Euler-Lagrange compact} we then find
\eqref{weak sol eq}.

   In order to prove the bound \eqref{eq:exbound}, we simply pass to the
limit as $\tau \to 0$ in \eqref{eq:prebound0} and obtain
$$ \|y \|_{L^\infty(0,T;W^{2,p}(U;\Rz^d))}^{   p} + \|\nabla \dot y^\top \nabla
y+ \nabla
y^\top \nabla \dot y  \|_{L^2(Q;\Rz^{d\times d})}^{   2} \leq c$$
independently of $\eps$. Following again \cite[Thm. 3.1]{HealeyKromer}
we have that $\det \nabla y \geq c>0$ independently of $\eps$. By
\cite{Neff} and \cite[Thm. 2.2]{Pompe} this ensures that
$$\| \nabla \dot y \|_{L^2(Q;\Rz^{d\times d})}^2\leq c \|\nabla \dot y^\top \nabla
y+ \nabla
y^\top \nabla \dot y  \|_{L^2(Q;\Rz^{d\times d})}^2\leq c$$
independently of $\eps$. Hence, \eqref{eq:exbound} follows by the
Poincar\'e inequality. \end{proof}

Before moving to the proof of Theorem \ref{Thm:existence sol} in the
diffused-interface case $\eps>0$, let us recall a well-posedness
result for the growth subproblem, see \cite[Thm.~3.15]{Mennucci}.

\begin{proposition}[Well-posedness of the growth
  problem] \label{prop:mennucci}
  Assume to
  be given $\haz \gamma \in C(\Rz^d)$ with $c_\gamma \leq
  \haz \gamma(\cdot) \leq C_\gamma$ for some $0<c_\gamma\leq C_\gamma$ and
  $\Omega_0\subset \Rz^d$ nonempty, open, and bounded. Then, there
  exists a unique viscosity solution to
  \begin{align}
    &\haz \gamma (x) |\nabla (-\theta)(x)|=1 \quad \text{in} \ \ \Rz^d
    \setminus \overline{\Omega_0},\label{eq:new1}\\
    &\theta=0 \quad \text{in} \  \ {\Omega_0}. \label{eq:new2} 
  \end{align}
 Moreover, $\theta \in C^{0,1}(\Rz^d)$ with  
  \begin{equation}
  \label{eq:lip}
  0<\frac{1}{C_\gamma} \leq |\nabla \theta(x)| \leq \frac{1}{c_\gamma} \ \
  \text{for a.e.} \ \ x\in \Rz^d,
\end{equation}
and we have that
\begin{equation}\label{bound theta with distance}
    \frac{\operatorname{dist}(x,{\Omega_0})}{C_\gamma}\leq \theta(x)
    \leq \frac{\operatorname{dist}(x,{\Omega_0})}{c_\gamma} \quad
    \forall x \in   \Rz^d  \setminus \overline{\Omega_0}.
  \end{equation}
\end{proposition}

We are now ready to prove   Theorem \ref{Thm:existence sol} in the diffused-interface case $\eps>0$. As announced, the proof hinges on an iterative construction. To
start with, let us remark that $y_0$ from \eqref{I1} is such that
$\nabla y_0$ is H\"older continuous. In particular,
the mapping $\tilde \gamma :    \overline{U}  \to (0,\infty)$ defined by
$$\tilde \gamma (   x ) := \gamma(y_0(   x ),\nabla y_0
(   x ))\quad \forall    x \in    \overline{U} $$
is H\"older continuous, as well.    Letting $\haz \gamma$ be any
continuous extension of $\tilde \gamma$ to $\Rz^d$ with $c_\gamma \leq
\haz \gamma(\cdot)\leq C_\gamma$, we can use
Proposition \ref{prop:mennucci} and find   $\theta_0
\in C(\overline U )$ solving  
\begin{align*}
    &\gamma(y_0(    x ),\nabla y_0(    x ))|\nabla (-\theta_0)(x)|=1\quad \text{in} \  U  \setminus\overline{\Omega_0},\\
    &   \theta_0 =0 \quad \text{in} \  \Omega_0
\end{align*}
   in the viscosity sense, with \eqref{eq:lip} and \eqref{bound
  theta  with distance} holding in $\overline U$.   Note that \eqref{bound theta with distance} in particular implies that
     $$\Omega^0(T)=\{x\in U \ | \ \theta_0(x)<T\} \subset
 \Omega_0 +B_{C_\gamma T}\stackrel{\eqref{hp dist Omega0 to partial B}}{\subset\subset} U.$$
By applying Proposition \ref{prop: existence y} for
$\theta=\theta_0$ we find $y^1
\in  L^{\infty}(0,T;W^{2,p}(U ;\R^d))\cap H^1(0,T;H^1(U
;\R^d))$.

This can be iterated as follows: For    all $k\geq 1$, given $y^k \in  L^{\infty}(0,T;W^{2,p}(U ;\R^d))\cap H^1(0,T;H^1(U
  ;\R^d)) $ we  define 
$\theta^k\in C(\overline U )$ to  be a viscosity
solution to 
  \begin{align*}
   &   \gamma(y^k(\theta^k (x)   {\wedge} T  ,x),\nabla
     y^k(\theta^k(x)   {\wedge} T  ,x))|\nabla (-\theta^k) (x)|=1\quad \text{in} \  U  \setminus\overline{\Omega_0},\\
      &   \theta^k   =0  \quad \text{in} \  \Omega_0
  \end{align*}
     with \eqref{eq:lip} and \eqref{bound
  theta  with distance} holding in $\overline U$.  
   The existence of such a viscosity solution follows again from
     Proposition \ref{prop:mennucci} as the mapping on $ \overline{U} $ defined   as
$$    x \mapsto \gamma(y^k(\theta^k (x)   {\wedge} T  ,x),\nabla
     y^k(\theta^k(x)   {\wedge} T  ,x))\quad \forall x \in \overline{U}$$
        may be extended to a continuous mapping $\haz \gamma$ on
     $\Rz^d$   with $c_\gamma \leq
\haz \gamma(\cdot)\leq C_\gamma$.   
        Note again that \eqref{bound theta with distance}   implies that
     \begin{equation}\label{eq:notouch}\Omega^k(   T   )
 :=\{x\in U \ | \ \theta^k(x)<   T  \} \subset
 \Omega_0 +B_{C_\gamma T}  \stackrel{\eqref{hp dist Omega0 to partial B}}{\subset\subset} U.
\end{equation}
Inclusion   \eqref{eq:notouch} in
 particular guarantees that the accreting phase defined by
 $\theta^k$ remains at positive distance from the boundary
 $\partial U$, independently of $\eps$ and $k$. 

 Given such $\theta^k$, we  define $y^{k+1} \in
   L^{\infty}(0,T;W^{2,p}(U ;\R^d))\cap H^1(0,T;H^1(U ;\R^d))
  $ by  Proposition \ref{prop: existence y}  applied  for
  $\theta=\theta^k$.  

   Bounds \eqref{eq:exbound} and    \eqref{eq:lip}     ensure that   
the sequence $(y^k, \theta^k)_{k\in \N}$
defined by this iterative  procedure  is 
(possibly not unique but nonetheless)
uniformly  bounded in 
$$\left(L^{\infty}(0,T;W^{2,p}(U ;\R^d))\cap H^1(0,T;H^1(U
  ;\R^d))\right) 
\times
  C^{0,1}(\overline U). $$ 
   As   
$(\theta^k)_{i\in\N}$    are    uniformly Lipschitz continuous, 
by the
Ascoli--Arzel\`a   and  the Banach--Alaoglu Theorems,    possibly
passing to not relabeled subsequences,    
   one can find   a pair  $(y,
\theta) 
$
such that
\begin{align}
  &y^k \stackrel{*}{\rightharpoonup} y \quad \text{weakly-$*$ in} \  L^{\infty}(0,T;W^{2,p}(U ;\R^d))\cap H^1(0,T;H^1(U  ;\R^d)),\label{eq:conv_y}\\
  &y^k \to y \quad \text{strongly in} \  C^{1,\alpha}(
    \overline{Q}  ;\R^d), \label{eq:conv_y2}\\
  &\theta^k \to \theta \quad \text{strongly in} \  C( \overline{U}) \label{eq:conv_theta}
\end{align}
 for some $\alpha \in (0,1)$ and $\theta$ fulfills    \eqref{eq:lip}
 and \eqref{bound theta  with distance} in $\overline U$.  
 As $(y^k)_{k\in \N}$    are   uniformly H\"older continuous and
$\gamma$ is Lipschitz continuous, by \eqref{bounds gamma} we have
\begin{align*}
  &|\gamma(y^k(\theta^k(x)   {\wedge} T  ,x),\nabla
  y^k(\theta^k(x)   {\wedge} T  ,x ) )-
  \gamma(y^j(\theta^j(x)   {\wedge} T  ,x),\nabla
    y^j(\theta^j(x)   {\wedge} T  ,x ))|\\
  &\quad \leq
  c|y^k(\theta^k(x)   {\wedge} T  ,x)-y^j(\theta^j(x)   {\wedge} T  ,x)|+c|\nabla
  y^k(\theta^k(x)   {\wedge} T  ,x)-\nabla y^j(\theta^j(x)   {\wedge} T  ,x)|\\
  &\quad \leq
     c|y^k(\theta^k(x)   {\wedge} T  ,x)-y^j(\theta^k(x)   {\wedge} T  ,x)|+c|\nabla
  y^k(\theta^k(x)   {\wedge} T  ,x)-\nabla y^j(\theta^k(x)   {\wedge} T  ,x)|\\
  &\qquad + c|y^j(\theta^k(x)   {\wedge} T  ,x)-y^j(\theta^j(x)   {\wedge} T  ,x)|+c|\nabla
    y^j(\theta^k(x)   {\wedge} T  ,x)-\nabla y^j(\theta^j(x)   {\wedge} T  ,x)|\\
  &\quad \leq c\| y^k - y^j\|_{C^1( \overline Q    ; \Rz^d   ) }
    + c \| \theta^k - \theta^j\|^\alpha_{C(\overline
    U)} \quad  \forall \ x \in \overline U.
\end{align*}
   Together with  
\eqref{eq:conv_y2}--\eqref{eq:conv_theta}, this proves that  $x \mapsto \gamma(y^k(\theta^k(x)   {\wedge} T  ,x),\nabla
  y^k(\theta^k(x)   {\wedge} T  ,x ) )$ converges to $x \mapsto \gamma(y(\theta(x)   {\wedge} T  ,x),\nabla
  y(\theta(x)   {\wedge} T  ,x ) )$ uniformly in $\overline U$. By  the stability
  of the eikonal equation with respect to the uniform convergence of
  the data,  see, e.g., \cite[Prop.~1.2]{Ishii},   $\theta$ satisfies
  \eqref{eq: eikonal system}--\eqref{eq: eikonal system2}  
  with   coefficient  $x \mapsto \gamma(y(\theta(x)   {\wedge} T  ,x),\nabla
  y(\theta(x)   {\wedge} T  ,x ) )$. Moreover, since    bound
  \eqref{eq:exbound} is 
  independent of $\theta$, following the argument of the proof
of Proposition \ref{prop: existence y}, we can pass to the limit in
the Euler--Lagrange equation \eqref{weak sol eq} and conclude the proof of Theorem \ref{Thm:existence sol} in the
case $\eps>0$.

\section{Proof of Theorem \ref{Thm:existence sol}:  sharp-interface case}\label{Sec: eps to zero}

The existence of weak/viscosity solutions in the
sharp-interface case $\eps=0$ is obtained by 
passing to the limit as $\eps\to 0$ in sequences of weak/viscosity
solutions $(y_\eps,\theta_\eps)$ of the diffused-interface
problem.

Notice at first that $\theta_\eps$ are uniformly Lipschitz
continuous, see \eqref{eq:lip}.    Bound
\eqref{eq:exbound} is  
independent of $\eps$ and    implies   that
there exist not relabeled subsequences such that
\begin{align}
  &y_\eps \stackrel{*}{\rightharpoonup} y \quad \text{weakly-$*$ in} \  L^{\infty}(0,T;W^{2,p}(U ;\R^d))\cap H^1(0,T;H^1(U  ;\R^d)),\label{y_eps conv}\\
  &y_\eps \to y \quad \text{strongly in} \  C^{1,\alpha}(\overline{Q};\R^d), \label{y_eps conv2}\\[1mm]
  &\theta_\eps \to \theta \quad \text{strongly in} \  C(\overline{U}) \label{theta_eps conv}
\end{align}
for some $\alpha\in (0,1)$.

   Let us now   prove that we can pass to the limit
  $\eps\rightarrow0$ in equation \eqref{weak sol eq}. The convergence
  of the loading    is straightforward.   Moreover, the level sets $\{ \theta(x)=t\}$ have Lebesgue measure zero by \eqref{eq:lip}.
  Hence, by the assumptions \eqref{h_eps properties} on $h_\eps$  and
  the uniform convergence \eqref{theta_eps conv} of
  $(\theta_\eps)_\eps$, we have    that   
  \begin{align*}
    h_\eps(\theta_\eps(x){-}t)\rightarrow h_0(\theta(x){-}t) \quad
    \text{ for a.e. } (t,x)\in  Q,
  \end{align*}
  and that $(t,x)\mapsto h_\eps(\theta_\eps(x){-}t)$ converges to
  $(t,x)\mapsto h_0(\theta(x){-}t)$ strongly in $L^2( Q )$.
On the other hand, by \eqref{W1}     and  
convergence    \eqref{y_eps conv2},   
   for all   $(t,x)\in Q $ and 
$i=   a,\, r,\,    J$,   we have that 
\begin{equation*}
  |\partial_F    V^i   (\nabla y_\eps)|  \leq c
  \quad\text{   and } \quad    \partial_F    V^i   (\nabla
  y_\eps)\rightarrow \partial_F    V^i   (\nabla y).
\end{equation*}
   Fix $z\in C^{\infty}([0,T]\times \overline{U };\R^d)$ with $
          z=0 $ on $\Sigma_D$. By Lebesgue's Dominated Convergence
          Theorem    we get  
\begin{align*}
  &\int_0^T\!\!\int_U \partial_F W_\eps(\theta_\eps{-}t,\nabla
    y_\eps)  {:}\nabla z   \,\dx\,\dt\\
   &\quad    =\int_0^T\!\!\int_U \Big(h_\eps(\theta{-}t)\partial_F   V^r
    (\nabla y_\eps) + (1-h_\eps(\theta{-}t))\partial_F   V^a
    (\nabla y_\eps) + \partial_F \J(\nabla y_\eps)\Big)   {:}\nabla z\,\dx\,\dt\\
  &\quad
    \rightarrow  \int_0^T\!\!\int_U  \Big(h_0(\theta{-}t)\partial_F  
    V^r  (\nabla y) +  (1{-}h_0(\theta{-}t))\partial_F   V^a
     (\nabla y) +    \partial_F \J(\nabla y)\Big)   {:}\nabla z  
    \,\dx\,\dt\\
  &\quad    =\int_0^T\!\!\int_U \partial_F W_0(\theta{-}t,\nabla
    y)   {:}\nabla z \,\dx\,\dt  
\end{align*}
Furthermore, by using convergence \eqref{y_eps conv}, we get  
\begin{align*}
 & \int_0^T\!\!\int_U \partial_{\dot{F}} R_\eps(\theta_\eps{-}t,\nabla
   y_\eps, \nabla \dot{y}_\eps){:}\nabla z \,\dx \,\dt \\
  &\quad \rightarrow \int_0^T\!\!\int_U\Big(
    h_0(\theta{-}t)\partial_{\dot{F}}R^r(\nabla y, \nabla
    \dot{y})+
    (1{-}h_0(\theta{-}t))\partial_{\dot{F}}R^{a}(\nabla y, \nabla
    \dot{y})\Big) {:}\nabla z \,\dx \,\dt\\
  & \quad    = \int_0^T\!\!\int_U \partial_{\dot{F}} R_0(\theta{-}t,\nabla
   y, \nabla \dot{y}){:}\nabla z \,\dx \,\dt  .
\end{align*}

   In order to prove the convergence of the second-order term, we
set $z_\eps=y-y_\eps $ and recall that
$z_\eps \to 0$ strongly in $L^\infty(0,T;H^1(U;\Rz^d))$ and $z_\eps
\stackrel{*}{\rightharpoonup} 0$ weakly-$*$ in  $L^\infty(0,T;W^{2,p}(U;\Rz^d))$
in order to obtain
\begin{align}
  &\limsup_{\tau\rightarrow 0} \int_0^T\!\!\int_U({ \rm D}H(\nabla^2 y)-{
    \rm D} H(\nabla^2 y_\eps )){\vdots}\nabla^2 z_\eps \,\dx\,\dt\nonumber\\
  &\quad=\limsup_{\eps\rightarrow 0}\Bigg(\int_0^T\!\!\int_U{ \rm D}
    H(\nabla^2 y ){\vdots}\nabla^2
    z_\eps-f(\theta{-}t){\cdot}z_\eps\,\dx\,\dt
  \nonumber\\
  &\qquad+\int_0^T\!\!\int_U\Big(\partial_F
    W_\eps(\theta{-}t, \nabla
    y_\eps) {:}\nabla z_\eps +\partial_{\dot
    F}R_\eps\left(\theta{-}t, \nabla
    y_\eps, \nabla
    \dot y_\eps\right) {:}\nabla z_\eps \Big)\,\dx \,\dt \nonumber
    \Bigg)=0
\end{align}
   Owing to    \eqref{H32} this proves that $\nabla^2 y_\eps \to \nabla^2 y$ strongly in
$L^p(Q;\Rz^{d\times d\times d})$. We hence have
that ${\rm D}H(\nabla^2 y_\eps) \to {\rm D}H(\nabla^2 y)$ strongly in
$L^{p'}(Q;\Rz^{d\times d\times d})$, as well, and we  can pass to the limit as $\eps \to 0$ in \eqref{weak sol
  eq}.

In order to conclude the proof, we are left to check that $\theta$ is
a viscosity solution to \eqref{eq:  eikonal system}. This however
readily follows as $x\mapsto \gamma(y_\eps(\theta_\eps(x)\wedge
T,x),\nabla y_\eps(\theta_\eps(x)\wedge T,x))$ converges to $x\mapsto \gamma(y(\theta(x)\wedge
T,x),\nabla y(\theta(x)\wedge T,x))$  uniformly and the eikonal problem
is stable under uniform convergence of the data
\cite[Prop.~1.2]{Ishii}.

Before closing this section, let us explicitly remark that indeed
Proposition \ref{prop:energy} actually holds in the case $\eps=0$, as
well. In order to check it, one would need a slightly different, and
indeed simpler, a-priori estimate on the time-discrete
solutions. Based on such result, one could argue as in Section
\ref{sec: Diffusive model} by the same iterative procedure
in order to obtain an alternative proof of Theorem \ref{Thm:existence sol} in the
sharp-interface case.  

\section*{Acknowledgments}

This research was funded in
whole or in part by the Austrian Science Fund (FWF) projects 10.55776/F65,  10.55776/I5149,
10.55776/P32788,     and   10.55776/I4354,
as well as by the OeAD-WTZ project CZ 09/2023. For
open-access purposes, the authors have applied a CC BY public copyright
license to any author-accepted manuscript version arising from this
submission. Part of this research was conducted during a visit to the
Mathematical Institute of Tohoku University, whose warm hospitality is gratefully acknowledged.

\end{document}